\documentclass[11pt]{article}

\usepackage{amsmath, amssymb}
\usepackage{enumerate}
\numberwithin{equation}{section}

\allowdisplaybreaks

\usepackage{graphicx}
\usepackage{color}
\newcommand{\overbar}[1]{\overline{\mkern-4mu#1\mkern-0.0mu}}

\newcommand{\p}{\partial}
\newcommand{\e}{\epsilon}

\newcommand{\nn}{\nonumber}
\newcommand{\N}{\mathbb{N}}

\newcommand{\Q}{\mathbb{Q}}

\newcommand{\C}{\mathbb{C}}
\newcommand{\bt}{{\bf t}}

\newcommand{\bft}{\mathbf{t}}

\newcommand{\beq}{\begin{equation}}
\newcommand{\eeq}{\end{equation}}

\newcommand{\M}{\overbar{\mathcal{M}}}

\renewcommand{\H}{\mathcal{H}}

\newtheorem{dfn}{Definition}[section]
\newtheorem{lem}[dfn]{Lemma}
\newtheorem{prp}[dfn]{Proposition}
\newtheorem{thm}[dfn]{Theorem}
\newtheorem{rmk}[dfn]{Remark}
\newtheorem{cor}[dfn]{Corollary}
\newtheorem{emp}[dfn]{Example}
\newtheorem{cnj}[dfn]{Conjecture}

\newenvironment{prf}{\noindent {\it Proof} \ }{\hfill $\Box$}
\newenvironment{prfn}[1]{\noindent {\it Proof of #1} \ }{\hfill $\Box$}

\DeclareMathOperator{\Res}{\mathrm{res}}

\begin{document}

\title{The Loop Equation for Special Cubic Hodge Integrals}
\author{
{Si-Qi Liu$^{*}$, \quad Di Yang$^{**}$, \quad Youjin Zhang$^{*}$, \quad Chunhui Zhou$^{**}$}\\
{\small ${}^{*}$ Department of Mathematical Sciences, Tsinghua University, Beijing 100084, P.R.~China}\\
{\small ${}^{**}$ School of Mathematical Sciences, University of Science and Technology of China,} \\
{\small Hefei 230026, P.R.~China}
}
\date{}
\maketitle
\begin{abstract}
As the first step of proving the Hodge-FVH correspondence recently proposed in~\cite{LZZ},
we derive the Virasoro constraints and the Dubrovin--Zhang loop equation for special cubic Hodge integrals. 
We show that this loop equation has a unique solution, and provide a new algorithm for the computation of these Hodge integrals. 
We also observe the gap phenomenon for certain special cubic Hodge free energies.
\end{abstract}

\tableofcontents

\section{Introduction}

Let $\M_{g,n}$ be the moduli space of stable algebraic curves of genus~$g$ with~$n$ distinct marked points, 
where $g$ and~$n$ are non-negative integers satisfying the stability condition $2g-2+n>0$. 
For $1\le k\le n$ and $0\le j\le g$,
denote by~$\psi_k$ the first Chern class of the $k$-th tautological line bundle $\mathbb{L}_k$ on $\overbar{\mathcal{M}}_{g,n}$, and by~$\lambda_j$ the $j$-th Chern class of the Hodge bundle $\mathbb{E}_{g,n}$ on~$\overbar{\mathcal{M}}_{g,n}$. 
The rational numbers defined by the formula
$$
\int_{\overbar{\mathcal{M}}_{g,n}}
\psi_1^{i_1}\cdots\psi_n^{i_n}\lambda_{j_1}\cdots\lambda_{j_m}
$$
are called the Hodge integrals. These numbers take zero value unless the degree-dimension counting matches, i.e.
\beq \label{dd} i_1+\dots +i_n+ j_1+\dots +j_m = 3g-3+n. \eeq 
Denote by $\mathcal{C}_g(z):= \sum_{j=0}^g \lambda_j z^j$ the Chern polynomial of~$\mathbb{E}_{g,n}$. 
We will be particularly interested in the following class of Hodge integrals
defined via the cubic products of Chern polynomials, called the cubic Hodge integrals:
\beq\label{cubHo}
\int_{\overbar{\mathcal{M}}_{g,n}}
\psi_1^{i_1}\cdots\psi_n^{i_n} \mathcal{C}_g(-p)\mathcal{C}_g(-q)\mathcal{C}_g(-r), 
\eeq
where $p,q,r$ are complex parameters. 
These integrals are called \emph{special} if $p,q,r$ satisfy the following
{\it local Calabi--Yau condition}:
\beq\label{calabi-yau}
pq+qr+rp=0.
\eeq
These Hodge integrals are important in the localization technique of computing Gromov--Witten invariants for toric three-folds \cite{GP, LLZ, OP}. Their significance was also manifested by the  
Gopakumar--Mari\~no--Vafa conjecture regarding the Chern--Simons/string duality \cite{GV,MV}.

Let $\H=\mathcal{H}(\bt; p,q,r; \e)$ be the cubic Hodge free energy defined by
\begin{align}
& \mathcal{H}(\bt; p,q,r; \e) := \sum_{g\geq0}\e^{2g-2}\mathcal{H}_g(\bt; p,q,r)  ,  \nn \\
& \mathcal{H}_g(\mathbf{t}; p,q,r) := \sum_{n\geq 0} \sum_{i_1,\dots,i_n\geq 0} \frac{t_{i_1}\cdots t_{i_n}}{n!}
 \int_{\overbar{\mathcal{M}}_{g,n}}  \psi_1^{i_1}\cdots\psi_n^{i_n} \mathcal{C}_g(-p)\mathcal{C}_g(-q)\mathcal{C}_g(-r) . \nn
\end{align}
Here $\mathcal{H}_g(\mathbf{t}; p,q,r)$ is called the genus~$g$ part of the free energy~$\H$. Then the exponential
\[e^{\mathcal{H}(\bt;p,q,r;\e)} =:Z_{\rm cubic}(\bt;p,q,r;\e)=:Z_{\rm cubic}\] is called the cubic Hodge partition function. 
Clearly, $\H_g(\bt; p,q,r)\in \C[p,q,r][[\bt]]$. 
The genus zero free energy~$\H_0(\bt)$ is actually independent of~$p,q,r$ and has the explicit expression
$$
\H_0(\bt) = \sum_{n\geq 3} \frac{1}{n(n-1)(n-2)} \sum_{i_1+\dots+i_n=n-3} \frac{t_{i_1}}{i_1!} \cdots \frac{t_{i_n}}{i_n!} .
$$

Define 
\beq
v(\bt):=\p_{t_0}^2 \H_0(\bt)=\sum_{n\geq 1} \frac1n  \sum_{i_1+\dots+i_n=n-1} \frac{t_{i_1}}{i_1!} \cdots \frac{t_{i_n}}{i_n!}.\label{vexpression}
\eeq
It satisfies the following Riemann hierarchy:
\beq\label{kdv-genus-zero}
\frac{\p v}{\p t_i} = \frac{v^i}{i!} \frac{\p v}{\p t_0} , \quad i\geq 0.
\eeq
More generally, if one defines $w=\e^2 \p_{t_0}^2 \H(\bt;p,q,r;\e)$, then $w$ satisfies an integrable 
hierarchy of Hamiltonian evolutionary PDEs \cite{BPS1,BPS2, DLYZ-1}, called the Hodge hierarchy for the special cubic Hodge integrals, 
which is a deformation of the Riemann hierarchy. The first member of this integrable hierarchy reads
$$
w_{t_1} = w w_{t_0} + \frac{\e^2}{12} \left(w_{t_0t_0t_0}  - (p+q+r) w_{t_0} w_{t_0t_0} \right) + O(\e^4). 
$$

The Hodge-FVH correspondence is given by the following conjecture~\cite{LZZ}.
\begin{cnj}
The Hodge hierarchy for the special cubic Hodge integrals is equivalent, under a certain Miura type transformation, 
to the {\it fractional Volterra hierarchy} (FVH). Furthermore, the special cubic Hodge partition function 
gives a tau function of the FVH. 
\end{cnj}

For the case with $p=q$, the 
 validity of the Hodge-FVH correspondence is 
 implied by the Hodge-GUE correspondence established in~\cite{DLYZ-2, DY}.
The goal of this and the subsequent paper~\cite{LYZZ}
is to prove the Hodge-FVH correspondence. 
In the present paper, we derive the Dubrovin--Zhang loop equation by studying the 
Virasoro constraints for the special cubic Hodge partition function. 
We show that this loop equation together with the genus zero free energy 
uniquely determines the partition function. 
In the subsequent paper, we will show that the fractional Volterra hierarchy admits the same 
Virasoro constraints and that there exists a particular tau function of this integrable hierarchy 
which is uniquely determined by the same loop equation, and we prove in this way the 
Hodge-FVH correspondence.

From now on, we assume that $p,q,r$ satisfy the local Calabi--Yau condition~\eqref{calabi-yau}. 
The case with $p,q,r\in \Q$ is called {\it rational}. Note that 
the Virasoro constraints for the special cubic Hodge partition function in the general case are quite complicated. 
However, in the rational case, we find explicit expressions 
of the Virasoro constraints which lead to the Dubrovin--Zhang loop 
equation. It turns out that the loop equation for the general case can be deduced from 
the one for the rational case.

Let us first consider the rational case. Due to the symmetry property of the cubic Hodge integrals with respect to $p,q,r$,  
and the homogeneity property (deduced from \eqref{dd})
\[
\mathcal{H}_g( \mathbf{t}; \lambda p, \lambda q, \lambda r)|_{t_i\mapsto t_i \lambda^{i-1}}  
= \lambda^{3g-3} \mathcal{H}_g( \mathbf{t};  p,  q,  r),
\]
we can assume that
\beq\label{rational-condition}
p=\frac{1}{K_1},\quad q=\frac{1}{K_2},
\quad r=-\frac{1}{h},
\eeq
where $K_1, K_2 \in \N$, $(K_1,K_2)=1$ and $h:=K_1+K_2$. 
We denote, for $\ell\geq 0$, 
\begin{align}
& b_{\alpha+h\ell} := \frac{\alpha}{K_1}+\ell, \quad c_{\alpha+h\ell} := \binom{b_{\alpha+h\ell} h}{b_{\alpha+h\ell} K_1} , \qquad \alpha=0,\dots,K_1-1,\label{bcdef1}\\ 
& b_{\alpha+h\ell} := \frac{-\alpha}{K_2}+\ell, \quad c_{\alpha+h\ell} := \binom{b_{\alpha+h\ell} h}{b_{\alpha+h\ell} K_2}, \qquad \alpha=-(K_2-1),\dots,-1,\label{bcdef2}
\end{align} 
and  
\[ \N_*=(\N-K_2) \backslash \left(\{0\} \cup (h\N-K_2)\right), \quad \mbox{where } a\N-K_2:=\{ak-K_2|k\in \N\}.\]
Define
\beq\label{Zxs}
Z(x,{\bf s};\e):=\exp\left(\frac{A(x,\tilde {\bf s})}{\e^2}\right) 
Z_{\rm cubic}\left(\mathbf{t}(x,{\bf s}); \frac{1}{K_1}, \frac{1}{K_2}, -\frac{1}{h};\e \right), 
\eeq
where ${\bf s}:=(s_k)_{k\in \N_*}$ is an infinite vector of indeterminates, 
$\tilde{s}_k=s_k-c_h^{-1}\delta_{k,h}$~($k\in \N_*$),
\beq\label{t-s}
t_i=t_i(x,{\bf s})=\sum_{k\in \N_*}  b_k^{i+1} c_k \tilde s_k+\delta_{i,1}+x \delta_{i,0},\quad i\geq 0, 
\eeq
and~$A$ is the quadratic series
\beq
A:=A(x,{\bf s})=\frac{1}{2}\sum_{k_1,k_2\in \N_*} \frac{b_{k_1} b_{k_2}}{b_{k_1} + b_{k_2}} c_{k_1} c_{k_2} s_{k_1} s_{k_2} + x \sum_{k\in \N_*} c_k s_k.
\eeq 
Note that for $g\geq0$, 
$\mathcal{H}_g(\mathbf{t}(x,\mathbf{s});\frac{1}{K_1}, \frac{1}{K_2}, -\frac{1}{h})$ 
is a well-defined formal power 
series in $\C[[x-1]][[\mathbf{s}]]$. 
Indeed, the coefficient of each monomial $(x-1)^{k_0}s_{k_1}\cdots s_{k_m}$ in 
$\mathcal{H}_g(\mathbf{t}(x,\mathbf{s});\frac{1}{K_1}, \frac{1}{K_2}, -\frac{1}{h})$ is a finite sum because of the dimension reason~\eqref{dd}.

\begin{rmk}
The change of the variables \eqref{t-s} is not invertible, but we do not lose any information by making the substitution \eqref{t-s} 
as we will see from Theorem \ref{refineloopmain} (cf. Proposition~\ref{uniquenessprop}) that 
 the special cubic Hodge integrals are uniquely determined by the loop equation \eqref{loop-eq-0}.
\end{rmk}

Denote $I=\{-(K_2-1),\dots,K_1-1\}$ and $I_*=I\backslash \{0\}$, and define a family of linear operators $L_m=L_m\left(\e^{-1}x,\e^{-1}\mathbf{s},\e\p/\p \mathbf{s}\right)$, $m\geq 0$ by
\begin{align}
L_0=&\sum_{k\in \N_*} b_k s_k \frac{\p}{\p s_k}+\frac{x^2}{2\e^2}+\frac{1}{24}\left(\frac1h-\frac1{K_1}-\frac1{K_2}\right), \label{virasoro-target0} \\
L_m=&\sum_{k\in \N_*} b_k s_k \frac{\p}{\p s_{k+h m}}+x\frac{\p}{\p s_{hm}} \nn\\
& +\frac{\e^2}{2} \sum_{\ell=1}^{m-1}\frac{\p^2}{\p s_{h \ell} \p s_{h(m-\ell)}} + \frac{\e^2}{2}  \sum_{\alpha,\beta \in I_*} \sum_{\ell=0}^{m-1}  G^{\alpha\beta}  \frac{\p^2}{\p s_{\alpha+h\ell} \p s_{\beta+h(m-1-\ell)}},\label{virasoro-target} 
\end{align}
where $\left(G^{\alpha\beta}\right)_{\alpha,\beta\in I}$ is a symmetric nondegenerate constant matrix defined by
\begin{align}\label{Galphabeta}
G^{\alpha\beta}=
\begin{cases}
\frac{K_1}{h}\delta^{\alpha+\beta,-K_2}, &\quad \alpha,\beta<0;\\
\quad 1,                                       &\quad \alpha=\beta=0;\\
\frac{K_2}{h}\delta^{\alpha+\beta, K_1}, &\quad \alpha,\beta>0;\\
\quad 0,                                       &\quad \text{elsewhere}.
\end{cases}
\end{align}
It is easy to check that the operators $L_m$ satisfy the following Virasoro commutation relations:
$$
\left[ L_m,L_n \right]= \left(m-n\right) L_{m+n}, \quad \forall\, m,n\geq 0.
$$
\begin{thm}\label{thm-virasoro}
For the rational numbers $p,q,r$ given by~\eqref{rational-condition},
the series $Z(x,\mathbf{s};\e)$ defined by~\eqref{Zxs} satisfies the following Virasoro constraints:
\beq\label{VirasoroZ}
L_m\left(\e^{-1}x,\e^{-1}\tilde{\mathbf{s}},\e\p/\p \mathbf{s}\right)
Z(x,\mathbf{s};\e)=0, \quad m\geq 0.
\eeq
\end{thm}

Using Theorem~\ref{thm-virasoro} and the technique developed in~\cite{DZ}, 
we derive in Section~\ref{sec-loop} the Dubrovin--Zhang loop equation for the special cubic Hodge 
free energies in the rational case; see Theorem~\ref{main}.

We proceed  to the general case. Denote 
\beq\label{sigma12}
\sigma_1 = -(p+q+r), \quad \sigma_3 = -2 \left(p^3+q^3+r^3\right).
\eeq
From the local Calabi--Yau condition~\eqref{calabi-yau}
and the fact that the integral in~\eqref{cubHo} is symmetric in $p,q,r$,
we know that 
\beq\label{polynomiality}
\H_g:=\H_g(\bt;p,q,r)\in \C[\sigma_1,\sigma_3][[\bt]], \quad g\ge 0.
\eeq
The following theorem is the main result of the present paper.
\begin{thm} \label{refineloopmain} 
The equation
\begin{align} 
&\sum_{i\geq0 }\left( \p^i \Theta+\sum_{j=1}^i \binom{i}{j} P_{j-1,i-j+1}\right) \frac{\p \Delta H }{\p z_i} \nn\\
= \; & \frac{\Theta^2}{16}-\left(\frac{1}{16}-\frac{\sigma_1}{24}\right)\Theta+ 
\e^2\sum_{i\geq0}\p^{i+2}\left(\frac{\Theta^2}{16}-\left(\frac{1}{16}-\frac{\sigma_1}{24}\right)\Theta\right)
\frac{\p \Delta H}{\p z_i}\nn\\
& + \frac{\e^2}{2}\sum_{i,j\geq0} P_{i+1,j+1} \left(\frac{\p^2\Delta H}{\p z_i \p z_j}
+\frac{\p \Delta H}{\p z_i}\frac{\p \Delta H}{\p z_j}\right),\label{loop-eq-0}
\end{align}
which is called the Dubrovin--Zhang loop equation for the special cubic Hodge integrals,  
has a unique solution of the form
\[\Delta H:=\sum_{g\geq 1} \e^{2g-2} H_g,\quad H_g:=H_g(z_0,\dots,z_{3g-2};\sigma_1,\sigma_3)\]
up to the addition of a constant to each $H_g$, $g\geq 1$.
These constants can be uniquely determined by the following conditions:
\begin{align}
& H_1= \frac1{24} \log z_1 + \frac{\sigma_1}{24} z_0, \quad 
\sum_{j=1}^{3g-2} j z_j \frac{\p H_g}{\p z_j} = (2g-2)H_g, \quad g\geq 2.  \label{eq2233} 
\end{align}
Here 
\[\p:=\sum_{k\geq 0} z_{k+1} \frac{\p}{\p z_k},\quad \Theta:= \frac1{1-e^{z_0}/\mu},\]
the coefficients $P_{i,j}$ are certain polynomials in
$\Theta, \sigma_1,\sigma_3, z_1,z_2,\dots$
whose explicit expressions are given in Section~\ref{section4}, and $\mu$ is an arbitrary parameter. Moreover, 
let $v(\bt)$ be defined in~\eqref{vexpression}, then the genus~$g$~($g\geq 1$) special cubic 
Hodge free energy has the expression
\beq\label{HgHg}
\mathcal{H}_g
=H_g\left(v(\bt), \frac{\p v(\bt)}{\p t_0},\cdots,\frac{\p^{3g-2} v(\bt)}{\p t_0^{3g-2}};\sigma_1,\sigma_3\right).
\eeq
\end{thm}

One can recursively solve the loop equation to obtain the free energies $H_g$, $g\geq 1$.
For example, the first two $H_g$ are given by
\begin{align}
H_1=&\frac{1}{24}\log z_1+\frac{\sigma_1}{24} z_0, \label{explicitH1} \\
H_2
=&\frac{1}{1152 }\frac{z_4}{z_1^2} -\frac{7 }{1920 }\frac{z_2 z_3}{z_1^3}+\frac{1}{360}\frac{z_2^3}{z_1^4}
+\frac{\sigma_1 }{480 }\frac{z_3}{z_1}-\frac{11\sigma_1 }{5760}\frac{z_2^2}{z_1^2}
+\frac{7\sigma_1^2 }{5760}z_2\nn\\
&+\left(\frac{\sigma_1^3}{17280}-\frac{\sigma_2}{34560}\right)z_1^2. \label{explicitH2}
\end{align}
These expressions agree with the results of~\cite{DLYZ-1}.

Besides the above theorems, we present in the next proposition  
some properties of the Taylor coefficients of $\H_g(\bt(x,{\bf s});\frac{1}{K_1}, \frac{1}{K_2}, -\frac{1}{h})$, $g\geq 1$ 
in the rational case, and we give in Theorem~\ref{gaptheorem} (see below) the {\it gap phenomenon} for  
the special cubic Hodge integrals in the rational case with the additional requirement that 
one of $p, q$ is equal to~1. This type of results was used 
by M.-X.~Huang, A.~Klemm and S.~Quackenbush to compute the Gromov--Witten invariants 
of the quintic Calabi--Yau three-fold up to genus~$51$~\cite{HKQ}.
\begin{prp}\label{gapcor}  The free energies 
$\mathcal{F}_g(x,{\bf s}):=\H_g(\bt(x,{\bf s});\frac{1}{K_1}, \frac{1}{K_2}, -\frac{1}{h})$, $g\geq 1$ with $\bt(x,{\bf s})$ defined in~\eqref{t-s} 
have the following property:
\begin{align}
& \mathcal{F}_1(x,{\bf s}) - \frac{\sigma_1-1}{24} \log x  \nn \\ 
= & \sum_{m\geq 1}\sum_{k_1,\dots,k_m\in \N_*} \frac{c_{1; \, k_1,\dots,k_m}(\sigma_1,\sigma_3)}{m!} x^{b_{k_1}+\dots+b_{k_m}-m} 
s_{k_1} \cdots s_{k_m},  \label{gap1}\\
& \mathcal{F}_g (x,{\bf s}) - \frac{R_g(\sigma_1,\sigma_3)}{x^{2g-2}} \nn\\ 
= & \sum_{m\geq 1}\sum_{k_1,\dots,k_m\in \N_*} \!\!\!\!\!\! \frac{c_{g; \, k_1,\dots,k_m}(\sigma_1,\sigma_3)}{m!} x^{b_{k_1}+\dots+b_{k_m}-m-(2g-2)} 
s_{k_1}\cdots s_{k_m}, \quad g\geq 2. \label{gap2}
\end{align}
Here $c_{g;k_1,\dots,k_m}(\sigma_1,\sigma_3)$, $g\geq 1$ are certain functions of $\sigma_1, \sigma_3$, 
and $R_g(\sigma_1,\sigma_3)$, $g\geq 2$ are certain polynomials of~$\sigma_1,\sigma_3$ 
satisfying the condition 
\begin{equation}
\deg R_g\leq 3g-3,\quad \textrm{with}\  \deg \sigma_1=1, \, \deg \sigma_3=3,\label{zh-10}
\end{equation} 
and the degree $(3g-3)$-part of~$R_g$ is given by
\beq \label{Faber}
   \frac{(-1)^g}{2(2g-2)!}  \frac{|B_{2g}||B_{2g-2}|}{2g (2g-2)} \left(\frac{\sigma_1^3}3-\frac{\sigma_3}{6}\right)^{g-1},
\eeq
where $B_j$ denote the Bernoulli numbers. 
\end{prp}

Note that in~\eqref{gap1} and~\eqref{gap2} the indeterminates of $R_g(\sigma_1,\sigma_3)$ take the values
\[\sigma_1=\frac1{K_1+K_2}-\frac1{K_1}-\frac1{K_2},\qquad
\sigma_3=\frac2{(K_1+K_2)^3}- \frac2{K_1^3} - \frac{2}{K_2^3}.\]
The first two $R_g$, $g\geq 2$ have the expressions
\begin{align*}
& R_2=-\frac{1}{1440}+\frac{13\sigma_1}{5760}-\frac{7 \sigma_1^2}{5760}+\frac{\sigma_1^3}{17280}-\frac{\sigma_3}{34560},\\
& R_3=\frac{1}{181440}-\frac{107 \sigma_1}{362880}+\frac{145\sigma_1^2}{290304}-\frac{961 \sigma_1^3}{4354560}
+\frac{31 \sigma_3}{2177280}+\frac{113 \sigma_1^4}{4354560}-\frac{113\sigma_1\sigma_3}{8709120}\\
& \qquad -\frac{\sigma_1^6}{13063680}+\frac{\sigma_1^3\sigma_3}{13063680}-\frac{\sigma_3^2}{52254720}.
\end{align*}

When one of $p,q$ is equal to~1, 
we have the following theorem.
\begin{thm}\label{gaptheorem}
For the rational case satisfying the condition that one of $p,q$ is equal to~1,
the free energies 
$\mathcal{F}_g(x,{\bf s})$, $g\geq 1$ 
 satisfy the following gap condition:
\begin{align}
& \mathcal{F}_1(x,{\bf s})|_{{\bf s_{\rm II}}= {\bf 0}} - \frac{\sigma_1-1}{24} \log x \; \in \; \C[[x,{\bf s_{\rm I}}]],  
\label{stronggap1}\\
& \mathcal{F}_g (x,{\bf s})|_{{\bf s_{\rm II}}= {\bf 0}} - \frac{R_g(\sigma_1,\sigma_3)}{x^{2g-2}} \; \in \; \C[[x,{\bf s_{\rm I}}]], 
\quad g\geq 2. \label{stronggap2}
\end{align} 
Here ${\bf s_{\rm I}}= (s_h, s_{2h}, s_{3h},\dots)$ and ${\bf s_{\rm II}}= (s_k)_{k\in \N_*\backslash h\mathbb{N}}$.
 \end{thm}

\noindent The proofs of Proposition~\ref{gapcor} and Theorem~\ref{gaptheorem} are given in~\cite{LYZZ}.
We hope that the results will be useful in the study of Gromov--Witten invariants for toric Calabi--Yau varieties.

\paragraph{Organization of the paper}
In Section~2, we derive the explicit expressions of Virasoro constraints for~$Z(x,{\bf s}; \e)$.  
In Section~3, we derive the Dubrovin--Zhang loop equation for the special cubic Hodge free energies in the rational case. 
In Section~\ref{section4}, we prove Theorem~\ref{refineloopmain}.

\paragraph{Acknowledgements} 
We would like to thank Boris Dubrovin, Shuai Guo, Yongbin Ruan, and Don Zagier for several very helpful discussions. 
This work is partially supported by NSFC No.\,11771238, No.\,11671371.
The work of S.-Q.~Liu is also supported by the National Science Fund for Distinguished Young Scholars No.\,11725104. Part of the work of C.~Zhou was done while he was a Ph.D. student of Tsinghua University; he acknowledges Tsinghua University for excellent working conditions and financial supports.

\section{Virasoro constraints}
In this section we first give two versions of Virasoro constraints for~$Z_{\rm cubic}$, 
and then we prove Theorem~\ref{thm-virasoro}.
Denote by~$Z_{\rm WK}$ the Witten--Kontsevich partition function
$$
Z_{\rm WK}(\mathbf{t};\e)=\exp
\left(\sum_{g\geq0}\e^{2g-2}\sum_{n\geq 0}
\frac{t_{i_1}\cdots t_{i_n}}{n!}
\int_{\overbar{\mathcal{M}}_{g,n}}\psi_1^{i_1}\cdots\psi_n^{i_n}\right).
$$
It is well known that $Z_{\rm WK}(\bft;\e)$ satisfies the following Virasoro constraints~\cite{DVV,Kontsevich,Witten}:
\beq\label{vir-kdv}
L_m^{\rm KdV}\left(\e^{-1}\tilde{\bft},\e\p/\p \bft\right) Z_{\rm WK}(\bft;\e)=0,\quad
\eeq
where $\tilde{t}_i=t_i-\delta_{i,1}$, and $L_m^{\rm KdV}:=L_m^{\rm KdV}\left(\e^{-1}\bft,\e\p/\p \bft\right)$ are linear operators given by
\begin{align}
& L_{-1}^{\rm KdV} :=\sum_{i\geq1}t_i\frac{\p}{\p t_{i-1}}+\frac{t_0^2}{2\e^2}, \nn \\
& L_0^{\rm KdV}  :=\sum_{i\geq0}\frac{2i+1}{2} t_i\frac{\p}{\p t_i}+\frac{1}{16}, \nn \\
& L_m^{\rm KdV} :=\sum_{i\geq0}\frac{(2i+2m+1)!!}{2^{m+1}(2i-1)!!} t_i\frac{\p}{\p t_{i+m}} +\frac{\e^2}{2}\sum_{i+j=m-1}\frac{(2i+1)!!(2j+1)!!}{2^{m+1}}\frac{\p^2}{\p t_i\p t_j}, \quad m\ge 1. \nn
\end{align}
These operators satisfy the Virasoro commutation relations
$$
\left[L_m^{\rm KdV},L_n^{\rm KdV}\right]=(m-n)L_{m+n}^{\rm KdV},\quad \forall\,m,n\geq -1.
$$

\begin{lem}[\cite{FP}] \label{lemmaFP}
The partition function~$Z_{\rm cubic}(\mathbf{t};p,q,r;\e)$ has the expression
\beq\label{Z-cubic}
Z_{\rm cubic}(\mathbf{t};p,q,r;\e)
=e^{G(\e^{-1}\tilde{\bft},\e\p/\p \bft)} Z_{\rm WK}(\mathbf{t};\e),
\eeq
where the operator~$G=G\left(\e^{-1}\bft,\e\p/\p \bft\right)$ is defined by 
$$
G\left(\e^{-1}\bft,\e\p/\p \bft\right) = - \sum_{i\geq1}\frac{B_{2i}}{2i(2i-1)}\left(p^{2i-1}+q^{2i-1}+r^{2i-1}\right)D_i,
$$ 
$B_{2i}$ are Bernoulli numbers,  and $D_i=D_i\left(\e^{-1}\bft,\e\p/\p \bft\right)$ are given by
$$
D_i := - \sum_{j\geq0}t_j\frac{\p}{\p t_{j+2i-1}}
+\frac{\e^2}{2}\sum_{j=0}^{2i-2}\frac{\p^2}{\p t_j \p t_{2i-2-j}}, \quad i\geq 1.
$$
\end{lem}

The lemma below follows from equation~\eqref{vir-kdv} and Lemma~\ref{lemmaFP}. 

\begin{lem}[\cite{Z}] \label{lemmazhou}
Define a sequence of operators $L_m^{\rm cubic}=L_m^{\rm cubic}\left(\e^{-1}\mathbf{t},\e \p / \p \mathbf{t}\right)$ by
$$
L_m^{\rm cubic}\left(\e^{-1}\mathbf{t},\e \p / \p \mathbf{t}\right) := e^G \circ L_m^{\rm KdV}\left(\e^{-1}\mathbf{t},\e \p / \p \mathbf{t}\right) \circ e^{-G}, \quad m\geq -1.
$$
Then $L_m^{\rm cubic}$ satisfy the Virasoro commutation relations
$$
\Bigl[L_m^{\rm cubic},L_n^{\rm cubic}\Bigr]=(m-n)L_{m+n}^{\rm cubic},\quad m,n\geq -1.
$$
Moreover, $Z_{\rm cubic}(\mathbf{t};p,q,r;\e)$ satisfies the equations
\beq\label{Virasorocubic}
L_m^{\rm cubic}\left(\e^{-1}\tilde{\mathbf{t}},\e\p / \p \mathbf{t}\right) Z_{\rm cubic}(\mathbf{t};p,q,r;\e)=0,
\quad m\geq -1.
\eeq
\end{lem}

We call~\eqref{Virasorocubic} the Virasoro constraints for~$Z_{\rm cubic}(\mathbf{t};p,q,r;\e)$.
Note that these Virasoro operators $L_m^{\rm cubic}$ are in general quite complicated, and it is difficult 
to use~\eqref{Virasorocubic} directly for the computation of the cubic Hodge integrals. For example, 
the operator $L_0^{\rm cubic}$ has the explicit expression
$$
L_0^{\rm cubic} = L_0^{\rm KdV} - \sum_{k\geq 1} \frac{B_{2k}}{2k}  \frac{\sigma_{2k-1}}{(2k-2)!} D_k,
$$
where $\sigma_{2k-1}=-(2k-2)! \left(p^{2k-1}+q^{2k-1}+r^{2k-1}\right)$, 
and the corresponding constraint given in~\eqref{Virasorocubic} has the expression
\begin{align}
& \left(\frac32\frac{\p}{\p t_1} + \sum_{k\geq 1} \frac{B_{2k}}{2k} \frac{\sigma_{2k-1}}{(2k-2)!}  \frac{\p}{\p t_{2k}} \right) 
Z_{\rm cubic} \nn\\ 
& = \left(\sum_{i\geq0}\frac{2i+1}{2} t_i\frac{\p}{\p t_{i}} +\frac1{16}  
+ \sum_{k\geq 1} \frac{B_{2k}}{2k}  \frac{\sigma_{2k-1}}{(2k-2)!}  
\left(\sum_{j\geq0}t_j\frac{\p}{\p t_{j+2k-1}} 
 - \sum_{j=0}^{2k-2}\frac{\p^2}{\p t_j \p t_{2k-2-j}}\right)\right)
Z_{\rm cubic}. \nn
\end{align}
Following~\cite{DLYZ-2} we consider the following linear combinations of Virasoro operators:
\beq\label{vir-til}
\widetilde{L}_m^{\rm cubic}=\widetilde{L}_m^{\rm cubic}\left(\e^{-1}\mathbf{t},\e \p / \p \mathbf{t}\right)
:=\sum_{k\geq-1}\frac{m^{k+1}}{(k+1)!}L_k^{\rm cubic},\quad m\geq 0.
\eeq
As it is shown in~\cite{DLYZ-2}, the operators $\widetilde{L}_m^{\rm cubic}$ also
satisfy the Virasoro commutation relations
$$
\left[\widetilde{L}_m^{\rm cubic},\widetilde{L}_n^{\rm cubic}\right]=(m-n)\widetilde{L}_{m+n}^{\rm cubic},
\quad \forall\, m,n\geq 0.
$$
From~\eqref{Virasorocubic} it follows that
\beq\label{Virasorocubic2nd}
\widetilde{L}_m^{\rm cubic}\left(\e^{-1}\tilde{\mathbf{t}},\e \p / \p \mathbf{t}\right) Z_{\rm cubic}(\mathbf{t};p,q,r;\e)=0,
\quad m\geq 0.
\eeq
\noindent We call~\eqref{Virasorocubic2nd} 
the second version of the Virasoro constraints for~$Z_{\rm cubic}(\mathbf{t};p,q,r;\e)$.

Let us proceed to derive the explicit expressions of~$\widetilde{L}_m^{\rm cubic}$. 
To this end, we are to use the Givental quantization technique~\cite{G} to simplify the computation. 
Note that we will use a slightly different convention of notations from Givental's,   
and we refer to Appendix~A of~\cite{DLYZ-2} for the details.
Let $(\mathcal{V},\omega)$ denote the Givental symplectic space, where $\mathcal{V}$ is the space of Laurent series over~$\C$,
and $\omega$ is the bilinear form on~$\mathcal{V}$ defined by
\beq
\omega(f,g):=-\Res_{z=\infty}f(-z)g(z)\frac{dz}{z^2},\quad \forall\, f,g\in\mathcal{V}.
\eeq
Denote by $q_i, p_i$~$(i\geq 0)$ the Darboux coordinates associated to $\omega$, and 
define a family of infinitesimal symplectic transformations $l_k$, $k\geq -1$ on~$\mathcal{V}$ by 
\beq
l_k:=(-1)^{k+1}z^{3/2}\p_z^{k+1} z^{-1/2},\quad k\geq -1.
\eeq
Then we have the following formulae~\cite{G}:
\beq\label{GiventalKdVlm}
L_k^{\rm KdV}
=\widehat{l_k}\big|_{q_i\mapsto t_i,\p_{q_i}\mapsto\p_{t_i}, \, i\geq 0} +\frac{\delta_{k,0}}{16}.
\eeq
Similarly, define a family of infinitesimal symplectic transformations $d_i:=z^{1-2i}$, $i\geq 1$ on~$\mathcal{V}$.
Then we have
\beq\label{diDi}
D_i \left(\e^{-1}\bft,\e\p/\p\bft\right) =\widehat{d_i}\big|_{q_j\mapsto t_j,\p_{q_j}\mapsto\p_{t_j}, \, j\ge0}.
\eeq
Note that the above results of this section hold true for general $p,q,r$. 
The local Calabi--Yau condition~\eqref{calabi-yau} for $p, q, r$ will be needed in what follows.

Following~\cite{DLYZ-2} we have the following lemma.
\begin{lem}\label{phi-expan}
Assume that $p,q,r$ satisfy the local Calabi--Yau condition~\eqref{calabi-yau}.
Define
\beq\label{Psidefi}
\Psi(z)=\left(p^{1/p}q^{1/q}r^{1/r}\right)^{-z}
\sqrt{\frac{\Gamma(1-z/p)\Gamma(1-z/q)\Gamma(1-z/r)}
{\Gamma(1+z/p)\Gamma(1+z/q)\Gamma(1+z/r)}} .
\eeq
Then the asymptotic expansion of~$\log\Psi(z)$ as $z\rightarrow\infty$ within a properly chosen sector 
is given by
\beq\label{Psiasymptotic}
\log \Psi(z) \sim \pm \frac{\pi i}4 -\sum_{i=1}^{\infty}\frac{B_{2i}}{2i(2i-1)} \left(p^{2i-1}+q^{2i-1}+r^{2i-1}\right)z^{1-2i} .
\eeq
\end{lem}
\begin{prf}
One can obtain formula~\eqref{Psiasymptotic} by using Stirling's formula~\cite{WW} and the definition~\eqref{Psidefi} 
after a careful calculation.
\end{prf}

We note that 
 the constants $\pm \pi i/4$ in the above asymptotic formula \eqref{Psiasymptotic} 
 do not affect the results of Givental's quantization. For simplicity, denote by $\log\Phi(z)$ the following formal power series of $z^{-1}$: 
\beq
 \log \Phi(z) := -\sum_{i=1}^{\infty}\frac{B_{2i}}{2i(2i-1)} \left(p^{2i-1}+q^{2i-1}+r^{2i-1}\right)z^{1-2i}.
\eeq

\begin{lem}\label{quantization} The operators~$\widetilde{L}_m^{\rm cubic}$, $m\geq 0$ have the expressions
\begin{align}
\widetilde{L}_m^{\rm cubic}\left(\e^{-1}\mathbf{t},\e\p/\p\mathbf{t}\right)
= \left(zV_m(z)e^{-m\p_z}\right)^{\wedge} \Big|_{q_i\mapsto t_i,\p_{q_i}\mapsto\p_{t_i},\,i\geq 0}+\frac{m}{16}-\frac{p+q+r}{24}, \label{quantizationgivental}
\end{align}
where $V_m(z)$, $m\geq 0$ are given by
\beq\label{v}
V_m(z)=\frac{\Phi(z)}{\Phi(z-m)}\sqrt{\frac{1}{1-m/z}}\,.
\eeq
\end{lem}
\begin{prf}
Denote  
\[\widehat{\Phi}:=\exp\left(\left(\log\Phi(z)\right)^{\wedge}\big|_{q_i\mapsto t_i,\p_{q_i}\mapsto\p_{t_i},\,i\geq 0}\right).\]
It follows from formula~\eqref{diDi} and 
Lemma~\ref{lemmazhou} that
\[L_m^{\rm cubic}\left(\e^{-1}\mathbf{t},\e\p/\p \mathbf{t}\right)
=\widehat{\Phi}\circ L_m^{\rm KdV}\left(\e^{-1}\mathbf{t},\e \p/\p \mathbf{t}\right)\circ \widehat{\Phi}^{-1}, \quad m\geq 0.\]
Then by using~\eqref{GiventalKdVlm} we obtain
\begin{align*}
L_m^{\rm cubic}\left(\e^{-1}\mathbf{t},\e\p/\p \mathbf{t}\right)
=&\left(\Phi\circ l_m \circ \Phi^{-1}\right)^{\wedge} \Big|_{q_i\mapsto t_i,\p_{q_i}\mapsto \p_{t_i}}+\frac{1}{16}\delta_{m,0}-\frac{p+q+r}{24}\delta_{m,-1}.
\end{align*}
Note that the term~$-\frac{p+q+r}{24}\delta_{m,-1}$ comes from the cocycle~\cite{G} in the quantization of~$\Phi\circ l_m \circ \Phi^{-1}$.
The quantization formula~\eqref{quantizationgivental} then follows from~\eqref{vir-til} as well as the identity 
\[z^{3/2}\Phi(z)\circ e^{-m\p_z}\circ \frac{1}{\sqrt{z}\Phi(z)}=z\frac{\Phi(z)}{\Phi(z-m)}\sqrt{\frac{1}{1-m/z}} \, e^{-m\p_z}.\]
The lemma is proved.
\end{prf}

We are now ready to prove Theorem~\ref{thm-virasoro}, 
where we recall that $p,q,r$ are assumed to be rational numbers given by~\eqref{rational-condition}.

\smallskip

\begin{prfn}{Theorem~\ref{thm-virasoro}}
In order to prove the validity of~\eqref{VirasoroZ}, it suffices to show that
the following identities hold true:
\beq\label{Lmconjugate}
L_m\left(\e^{-1}x,\e^{-1}\tilde{\mathbf{s}},\e \p / \p \mathbf{s}\right)
=K^m e^{\frac{A(x,\tilde{{\bf s}})}{\e^2}}\circ\widetilde{L}_m^{\rm cubic}\left(\e^{-1}\tilde{\mathbf{t}},\e\p/\p\mathbf{t}\right)\circ e^{-\frac{A(x,\tilde {\bf s})}{\e^2}}, \quad m\ge 0,
\eeq
where 
\begin{equation}\label{zh-9}
K=h^h K_1^{-K_1} K_2^{-K_2},
\end{equation}
 and
\beq\label{ttttt}
\tilde t_i=\tilde  t_i(x,{\bf s})=\sum_{k\in \N_*}  b_k^{i+1} c_k \tilde s_k+x \delta_{i,0},
\quad i\ge 0.
\eeq
For simplicity we denote 
\[\bar{s}_k:=c_k s_k,\quad k\in \N_*.\]
Then the above substitution of the independent variables~\eqref{ttttt} reads
\beq \label{t-s-0}
\tilde t_i (x,{\bf s}) =  \sum_{k\in \N_*}b_k^{i+1}\tilde{\bar{s}}_k+x \delta_{i,0}, \quad \tilde{\bar{s}}_k = \bar{s}_k - \delta_{k,h}. 
\eeq
Since both sides of~\eqref{Lmconjugate} satisfy the Virasoro commutation relations, 
we only need to prove~\eqref{Lmconjugate} for $m=0,1,2$. 

\begin{lem}
The operator $\widetilde{L}_0^{\rm cubic}$ satisfies the following relation:
\beq
e^{\frac{A(x,\tilde{{\bf s}})}{\e^2}}\circ \widetilde{L}_0^{\rm cubic} \left(\e^{-1}\tilde{\mathbf{t}},\e\p/\p\mathbf{t}\right) \circ e^{-\frac{A(x,\tilde{{\bf s}})}{\e^2}}
=\sum_{k\in \N_*}b_k \tilde{\bar{s}}_k \frac{\p}{\p \bar{s}_k}+\frac{x^2}{2\e^2}+\frac{\sigma_1}{24}.
\eeq
\end{lem}
\begin{prf}
By using Lemma~\ref{quantization} we have
\begin{align*}
\widetilde{L}_0^{\rm cubic} \left(\e^{-1}\mathbf{t},\e\p/\p\mathbf{t}\right)
=\widehat{z} \mid_{q_i\mapsto t_i,\p_{q_i}\mapsto\p_{t_i}}+\frac{\sigma_1}{24}
=\sum_{i\geq1}t_i\frac{\p}{\p t_{i-1}}+\frac{t_0^2}{2\e^2}+\frac{\sigma_1}{24}\,.
\end{align*}
Under the substitution~\eqref{t-s-0}, we arrive at
\begin{align*}
&\widetilde{L}_0^{\rm cubic}\left(\e^{-1}\tilde{\mathbf{t}},\e\p/\p\mathbf{t}\right)\nn\\
=&\sum_{i\geq 1}\left(\sum_{k\in \N_*}b_k^{i+1}\tilde{\bar{s}}_k\right)\frac{\p}{\p t_{i-1}}
+\frac{1}{2\e^2}\left(\sum_{k\in \N_*} b_k\tilde{\bar{s}}_k +x\right)^2
+\frac{\sigma_1}{24}\\
=&\sum_{k\in \N_*}b_k \tilde{\bar{s}}_k \frac{\p}{\p \bar{s}_k}
+\frac{1}{2\e^2}\left(\sum_{k\in \N_*} b_k\tilde{\bar{s}}_k +x\right)^2
+\frac{\sigma_1}{24}.
\end{align*}
The lemma is proved by applying the conjugation by~$e^{\frac{A(x,\tilde{{\bf s}})}{\e^2}}$ to the above equality.
\end{prf}

\begin{lem} The operator $\widetilde{L}_1^{\rm cubic}$ satisfies the following relation:
\begin{align}
& e^{\frac{A(x,\tilde{{\bf s}})}{\e^2}}\circ 
\widetilde{L}_1^{\rm cubic}\left(\e^{-1}\tilde{\mathbf{t}},\e\p/\p\mathbf{t}\right)\circ e^{-\frac{A(x,\tilde{{\bf s}})}{\e^2}}\nn\\
=&\sum_{k\in \N_*}b_k V_1(-b_k) \tilde{\bar{s}}_k \frac{\p}{\p \bar{s}_{k+h}}
+x V_1(0) \frac{\p}{\p \bar{s}_h}
+\frac{\e^2}{2}\sum_{\alpha=1}^{K_1-1}
\frac{\Res_{z=b_{\alpha}}V_1(z)}{b_{\alpha}}
\frac{\p^2}{\p \bar{s}_{\alpha} \p \bar{s}_{K_1-\alpha}}\nn\\
& +\frac{\e^2}{2}\sum_{\alpha=-(K_2-1)}^{-1}
\frac{\Res_{z=b_{\alpha}}V_1(z)}{b_{\alpha}}
\frac{\p^2}{\p \bar{s}_{\alpha} \p \bar{s}_{-\alpha-K_2}}, \label{compare1}
\end{align}
where $V_1(z)$ is defined in~\eqref{v}.
\end{lem}

\begin{prf}
Since $p=1/K_1,q=1/K_2,r=-1/h$, we find by using the definition~\eqref{v} and 
Lemma~\ref{phi-expan} that $V_1(z)$ can be given by the asymptotic series (as $z\to\infty$) of a rational function, namely, 
$$
\frac{\prod_{i=1}^{K_1+K_2}\left(z-\frac{i}{K_1+K_2}\right)}
{\prod_{i=1}^{K_1}\left(z-\frac{i}{K_1}\right) \prod_{i=1}^{K_2}\left(z-\frac{i}{K_2}\right)} \; \sim \; V_1(z)\,.
$$
For simplicity the above rational function will still be denoted by~$V_1(z)$ as this will not introduce ambiguity. 
(Without further mentioning, similar notations will also be used for $V_m(z)$ 
as they are also asymptotic series of rational functions.)
Note that the rational function~$V_1(z)$ has simple poles at $b_{-(K_2-1)},\dots,b_{-1}$, $b_1,\dots,b_{K_1-1}$, and~$1$, 
therefore, 
\[V_1(z)=1 +\frac{\Res_{y=1}V_1(y)}{z-1} +\sum_{\alpha\in I_*} \frac{\Res_{y=b_{\alpha}}V_1(y)}{z-b_{\alpha}}.\]
Hence we obtain
$$
V_1(z) =
1+\sum_{n\geq1}z^{-n} \left(\Res_{y=1}V_1(y)  + \sum_{\alpha\in I_*} b_{\alpha}^{n-1}\Res_{y=b_{\alpha}}V_1(y) 
\right).
$$
Denote 
$a_{\alpha}:=b_{\alpha}^{-1} \Res_{z=b_{\alpha}}V_1(z)$, $\alpha\in I_*$.
By taking $m=1$ in~\eqref{quantizationgivental} we obtain
\begin{align*}
&\widetilde{L}_1^{\rm cubic}  
= \left(zV_1(z)e^{-\p_z}\right)^{\wedge}\Big|_{q_i\mapsto t_i,\,\p_{q_i}\mapsto \p_{t_i}}
+\frac{1}{16}+\frac{\sigma_1}{24}\\
=&\sum_{i\geq0}\left(\sum_{j=0}^{i+1} t_j\binom{i+1}{j}
+\sum_{j=0}^i t_j\sum_{n=1}^{i+1-j}(-1)^n\binom{i+1}{j+n} \left(\Res_{y=1}V_1(y) + \sum_{\alpha\in I_*} a_{\alpha}b_{\alpha}^n\right)\right)
\frac{\p}{\p t_i}\\
  &+\frac{\e^2}{2}\sum_{i,j\geq0}\left(\sum_{n=0}^{i+1} (-1)^n \binom{i+1}{n} \left(\Res_{y=1}V_1(y)+ \sum_{\alpha\in I_*} a_{\alpha} b_{\alpha}^{n+j+1}\right)\right)
\frac{\p^2}{\p t_i \p t_j}
+\frac{t_0^2}{2\e^2}+\frac{1}{16}+\frac{\sigma_1}{24}\,.
\end{align*}
Now by performing the variables substitution~\eqref{t-s-0} we arrive at
\begin{align*}
&\widetilde{L}_1^{\rm cubic}\left(\e^{-1}\tilde{\mathbf{t}},\e\p/\p\mathbf{t}\right) \\
=& \sum_{k\in \N_*}b_k V_1(-b_k) \tilde{\bar{s}}_k \frac{\p}{\p \bar{s}_{k+h}}
+x V_1(0) \frac{\p}{\p \bar{s}_h}  
+\frac{1}{2\e^2}\left(\sum_{k\in \N_*} b_k \tilde{\bar{s}}_k+x\right)^2 +\frac{1}{16}+\frac{\sigma_1}{24}\\
&+\frac{\e^2}{2}\sum_{\alpha=1}^{K_1-1}a_{\alpha} \frac{\p^2}{\p \bar{s}_{\alpha} \p \bar{s}_{K_1-\alpha}}
+\frac{\e^2}{2}\sum_{\alpha=-(K_2-1)}^{-1} a_{\alpha} \frac{\p^2}{\p \bar{s}_{\alpha} \p \bar{s}_{-\alpha-K_2}}\\
&+\sum_{\alpha=1}^{K_1-1}a_{\alpha}\frac{\p A(x,\tilde {\bf s})}{\p \bar{s}_{\alpha}}\frac{\p}{\p \bar{s}_{K_1-\alpha}}
 +\sum_{\alpha=-(K_2-1)}^{-1}a_{\alpha}\frac{\p A(x,\tilde {\bf s})}{\p \bar{s}_{\alpha}}\frac{\p}{\p \bar{s}_{-\alpha-K_2}} .
 \end{align*}
Here we used the identity
\[\sum_{j=0}^i\sum_{k=1}^{i+1-j}\binom{i+1}{j+k}\lambda^j \mu^k
=\mu\frac{(\lambda+1)^{i+1}-(\mu+1)^{i+1}}{\lambda-\mu}.\]
Then dressing the operator $\widetilde{L}_1^{\rm cubic}\left(\e^{-1}\tilde{\mathbf{t}},\e\p/\p\mathbf{t}\right)$ 
by~$e^{\frac{A(x,\tilde {\bf s})}{\e^2}}$ we obtain~\eqref{compare1} after a long but straightforward computation.
The lemma is proved.
\end{prf}

We can prove in a similar way the following lemma.
\begin{lem} 
The operator $\widetilde{L}_2^{\rm cubic}$ satisfies the relation
\begin{align}
& e^{\frac{A(x,\tilde{{\bf s}})}{\e^2}}\circ 
\widetilde{L}_2^{\rm cubic}\left(\e^{-1}\tilde{\mathbf{t}},\e\p/\p\mathbf{t}\right)\circ e^{-\frac{A(x,\tilde{{\bf s}})}{\e^2}} \nn \\
=&\sum_{k\in \N_*}b_k V_2(-b_k) \tilde{\bar{s}}_k \frac{\p}{\p \bar{s}_{k+2h}}
+x V_2(0) \frac{\p}{\p \bar{s}_{2h}} +\frac{\e^2}{2}b_h^{-1}\Res_{z=b_h}V_2(z)\frac{\p^2}{\p \bar{s}_h \p \bar{s}_h} \nn \\
&  +\frac{\e^2}{2}\sum_{\alpha=1}^{K_1-1}\sum_{\ell=0}^1 b_{\alpha}^{-1} \Res_{z=b_{\alpha}}V_2(z) 
\frac{\p^2}{\p \bar{s}_{\alpha+h \ell} \p \bar{s}_{K_1-\alpha+h(1-\ell)}} \nn \\
&  +\frac{\e^2}{2}\sum_{\alpha=-(K_2-1)}^{-1}\sum_{\ell=0}^1 b_{\alpha}^{-1} \Res_{z=b_{\alpha}}V_2(z) 
\frac{\p^2}{\p \bar{s}_{\alpha+h \ell} \p \bar{s}_{-\alpha-K_2+h(1-\ell)}}. \label{compare2} 
\end{align}
\end{lem}

By using Lemma~\ref{phi-expan} we can also obtain the following lemma. 

\begin{lem}\label{const}
The numbers $c_k$, $k\in \N_*$ defined in \eqref{bcdef1}--\eqref{bcdef2} have the following properties:
\begin{enumerate}
\item [{\rm (i)}]
For $k\in \N_*$ and $\ell\geq 1$, $c_{k+h \ell}/c_k=K^{\ell} V_{\ell}(-b_k)$.
\item [{\rm (ii)}]
For $\ell\geq 1$, $c_{h \ell}=K^{\ell} V_{\ell}(0)$.
\item [{\rm (iii)}]  For $m,n\geq0$, 
\begin{align}
& c_{\alpha+h m} \, c_{K_1-\alpha+ h n} = \frac{h}{K_2} \frac{K^{m+n+1}}{b_{\alpha+h m}}  \Res_{z=b_{\alpha+h m}}V_{m+n+1}(z), \quad \alpha=1,\dots,K_1-1,  \nn\\
& c_{\alpha+h m} \, c_{-\alpha-K_2+ h n} = \frac{h}{K_1} \frac{K^{m+n+1}}{b_{\alpha+h m}}  \Res_{z=b_{\alpha+h m}}V_{m+n+1}(z), \quad \alpha=-(K_2-1),\dots,-1,\nn \\
& c_{h (m+1)} \, c_{h (n+1)}=\frac{K^{m+n+2}}{b_{h (m+1)}}\Res_{z=h (m+1)}V_{m+n+2}(z).   \nn
\end{align}
\end{enumerate}
\end{lem}

By comparing \eqref{compare1}--\eqref{compare2} with the definition of the left-hand side of~\eqref{Lmconjugate} and by using Lemma~\ref{const} 
we find that~\eqref{Lmconjugate} is true for $m=1, 2$. Thus Theorem~\ref{thm-virasoro} is proved.
\end{prfn}

\section{Loop equation: the rational case} \label{sec-loop}
In this section, we derive the Dubrovin--Zhang loop equation for the special cubic Hodge integrals in the rational case, namely, we take 
\[p=1/K_1,\quad q=1/K_2,\quad r=-1/h,\] 
where $K_1$ and $K_2$ are coprime positive integers, and $h=K_1+K_2$.

Introduce a generating series for the operators~$L_m$, $m\geq0$ defined in \eqref{virasoro-target0}--\eqref{virasoro-target} as follows:
\beq\label{Ltwo}
L(\lambda)=\sum_{m\geq 0}\frac{L_m(\e^{-1}x,\e^{-1}\tilde{\mathbf{s}},\e\p/\p \mathbf{s})}{\lambda^{m+1}}.
\eeq
\begin{lem}\label{virasoro-total}
The generating series~$L(\lambda)$ can be represented as
\beq\label{identityll}
L(\lambda)=\biggl[J_1(\lambda)J_2(\lambda)+\frac{1}{2}J_2(\lambda)^2\biggr]_{{\rm reg},-}
+\frac{1}{\lambda}\left(\frac{x}{2\e^2}+\frac{\sigma_1}{24}\right),
\eeq
where the operators $J_1(\lambda)$ and $J_2(\lambda)$ are defined by
\begin{align}
& J_1(\lambda) := \e^{-1}\frac{x}{\sqrt{\lambda}} + \e^{-1} \sum_{k \in \N_*}b_k \lambda^{b_k-\frac{1}{2}} \tilde{s}_k, \nn\\
& J_2(\lambda):=  \sum_{\ell\ge 1}\frac{\e}{\lambda^{b_{h\ell}+1/2}}\frac{\p}{\p s_{h \ell}} 
+ \sqrt{\frac{K_2}{h}}\sum_{\alpha=1}^{K_1-1} \sum_{\ell\geq0}\frac{\e}{\lambda^{b_{\alpha+h\ell}+\frac{1}{2}}}\frac{\p}{\p s_{\alpha+h\ell}} \nn\\
&\qquad \qquad \qquad + \sqrt{\frac{K_1}{h}}\sum_{\alpha=-(K_2-1)}^{-1}\sum_{\ell\geq0}\frac{\e}{\lambda^{b_{\alpha+h\ell}+\frac{1}{2}}}\frac{\p}{\p s_{\alpha+h\ell}}, \nn
\end{align}
and $[\,\cdot\,]_{{\rm reg},-}$ means to take terms of the series with negative integer powers of~$\lambda$.
\end{lem}
\begin{prf}
By comparing the coefficients of~$\lambda^{-(m+1)}$ for $m\geq 0$ on both
 sides of~\eqref{identityll} with $L_m$ given by \eqref{virasoro-target0}--\eqref{virasoro-target} we obtain the validity of the identity~\eqref{identityll}.
 The lemma is proved.
\end{prf}

Theorem~\ref{thm-virasoro} implies that $L(\lambda)Z(x,\mathbf{s};\e)=0$ holds true {\it identically} in~$\lambda$.
From the definition of~$Z(x,\mathbf{s};\e)$ given by~\eqref{Zxs}
we know that its logarithm admits the genus expansion
\[\log Z(x,\mathbf{s};\e)=\sum_{g\geq0}\e^{2g-2}\mathcal{F}_g(x,\mathbf{s}).\]
Here we denote, as we do in Proposition~\ref{gapcor},
\begin{equation}\label{zh-1}
\mathcal{F}_g(x, {\bf s}) :=
\mathcal{H}_g\left(\bt(x,{\bf s});\frac{1}{K_1}, \frac{1}{K_2}, -\frac{1}{h}\right)+A(x,\tilde{\bf s}) \, \delta_{g,0},\quad \textrm{with}\ \tilde s_k=s_k-c_h^{-1}\delta_{k,h}.
\end{equation}
In what follows we will also use the notations
\beq 
u(x,\mathbf{s}):=\p_x^2\mathcal{F}_0(x,\mathbf{s}) 
\eeq
and 
\begin{equation}\label{zh-5}
\Delta \mathcal{F}=\sum_{g\geq1}\e^{2g-2}\mathcal{F}_g(x,\mathbf{s}).
\end{equation}
\begin{lem}\label{prp-1}
The function $u(x, \mathbf{s})$ satisfies the following equations
\begin{align}
&\p_x^i u(x,\mathbf{s})=\p^i_{t_0}v(\mathbf{t}(x,\mathbf{s})), \quad i\ge 0.\\
&\frac{\p u}{\p s_k}=c_k \p_x e^{b_k u},
\quad \frac{\p^2\mathcal{F}_0}{\p x \p s_k}=c_k e^{b_k u},
\quad k \in \N_*.\label{dispersionlesshierarchy}
\end{align}
\end{lem}
\begin{prf}
By using~\eqref{t-s} we obtain
$$
\p_x=\p_{t_0},\quad 
\frac{\p}{\p s_k}=c_k \sum_{i\geq0}b_k^{i+1}\frac{\p}{\p t_i},\quad k\in\N_*.
$$
The lemma then follows from the Riemann hierarchy \eqref{kdv-genus-zero}.
\end{prf}

\begin{lem}[\cite{DLYZ-1}]\label{jetvariablesrep}
For any given $g\ge 1$, there exists a polynomial $H_g(z_0,\dots,z_{3g-2};p,q,r)$ of $p, q, r, z_2, z_3,\dots$ with coefficients depending smoothly on $z_0$ and rationally on~$z_1$, such that 
\beq
\mathcal{H}_g(\mathbf{t};p,q,r)=
H_g\left(v(\bt), \frac{\p v(\bt)}{\p t_0},\dots,\frac{\p^{3g-2} v(\bt)}{\p t_0^{3g-2}};p,q,r\right). \label{HHgHg}
\eeq
\end{lem}

It should be noted that Lemma~\ref{jetvariablesrep} does not require $p,q,r$ satisfy the local Calabi--Yau condition.

Let us denote 
\[u^{(i)}=\p_x^i u(x,{\bf s}),\quad i\geq 0.\]
Then from~Lemmas~\ref{prp-1}--\ref{jetvariablesrep} we know that 
\begin{equation}\label{zh-6}
\mathcal{F}_g(x,{\bf s})=H_g\left(u^{(0)}, u^{(1)},\dots,u^{(3g-2)}; p,q,r\right),\quad
g\geq 1, 
\end{equation}
where $p=1/K_1$, $q=1/K_2$, $r=-1/(K_1+K_2)$. 
Introduce a derivation~$D(\lambda)$ on $\C[[x-1,\mathbf{s}]]$ by
$$
D(\lambda)=
-\left[J_1(\lambda)J_2(\lambda)\right]_{{\rm reg},-}- \e^{-2} \left[J_2(\lambda)(\mathcal{F}_0)J_2(\lambda)\right]_{{\rm reg},-},
$$
where $J_1(\lambda), J_2(\lambda)$ are defined in Lemma \ref{virasoro-total}.

\begin{lem}\label{lemma35}
The series $D(\lambda)(\Delta\mathcal{F})$ has the following expression:
\begin{align}
&D(\lambda)(\Delta\mathcal{F}) \nn\\
=&\frac{\sigma_1}{24\lambda} + \biggl[
 \frac{1}{2\e^2} J_2(\lambda)^2(\mathcal{F}_0) +
\frac{1}{2}\sum_{i\geq 0}J_2(\lambda)^2(u^{(i)})\frac{\p \Delta\mathcal{F}}{\p u^{(i)}}  \nn \\
& +\frac{1}{2}\sum_{i,j\geq0}J_2(\lambda)(u^{(i)}) J_2(\lambda)(u^{(j)})
\left(\frac{\p^2 \Delta\mathcal{F}}{\p u^{(i)}\p u^{(j)}}
+\frac{\p \Delta\mathcal{F}}{\p u^{(i)}}\frac{\p \Delta\mathcal{F}}{\p u^{(j)}}\right)\biggr]_{{\rm reg},-}.  \label{loop-initial}
\end{align}
\end{lem}
\begin{prf}
From Lemma~\ref{virasoro-total} we know that the Virasoro constraints 
$L(\lambda)Z(x,\mathbf{s};\e)=0$ can be represented as
\begin{align}
&\biggl[J_1(\lambda)J_2(\lambda)(\mathcal{F}_0)
+\frac{1}{2 \e^2} (J_2(\lambda)(\mathcal{F}_0))^2\biggr]_{{\rm reg},-}
+\frac{x^2}{2\lambda}=0,\label{symmetry-0}\\
&\biggl[J_1(\lambda)J_2(\lambda)(\Delta\mathcal{F}) + \frac1{\e^2}J_2(\lambda)(\mathcal{F}_0)J_2(\lambda)(\Delta\mathcal{F})\nn\\
& \quad +\frac{1}{2}\left(J_2(\lambda)^2(\Delta\mathcal{F})
+(J_2(\lambda)(\Delta\mathcal{F}))^2\right)
+\frac{1}{2\e^2}J_2(\lambda)^2(\mathcal{F}_0)\biggr]_{{\rm reg},-}
+\frac{\sigma_1}{24\lambda}=0\label{symmetry-1}.
\end{align}
It is easy to see that~\eqref{symmetry-1} can be rewritten in the form
\beq
D(\lambda)(\Delta\mathcal{F})
=\frac{\sigma_1}{24\lambda} + \left[\frac{1}{2} \left(J_2(\lambda)^2(\Delta\mathcal{F})
+(J_2(\lambda)(\Delta\mathcal{F}))^2\right)
+\frac{1}{2\e^2}J_2(\lambda)^2(\mathcal{F}_0)\right]_{{\rm reg},-}.
\eeq
Since $J_2(\lambda)$ is a derivation, 
we obtain~\eqref{loop-initial}. The lemma is proved.
\end{prf}

Let us proceed to simplify the equation~\eqref{loop-initial}.
Introduce the following Puiseux series:
\begin{align} \label{eeplus}
E_{\alpha}=
\Bigg\{
\begin{array}{ll}
\sum_{n\geq0} c_{\alpha+hn} \zeta^{n+\frac{\alpha}{K_1}},
& \alpha=0,\dots,K_1-1,  \\
&\\
\sum_{n\geq0} c_{\alpha+hn} \zeta^{n-\frac{\alpha}{K_2}},
&  \alpha=-(K_2-1),\dots,-1,
\end{array}
\end{align}
where $\zeta:=e^{u(x,{\bf s})}/\lambda$. 
Define 
\beq\label{zh-2}
B_{i,j} := \sum_{\alpha,\beta\in I} \p_x^i \left(E_{\alpha}\right) G^{\alpha\beta}\p_x^j \left(E_{\beta}\right),
\eeq
where $G^{\alpha\beta}$ are constants given by~\eqref{Galphabeta}. 
Note that $B_{i,j}\in \C[u^{(1)},u^{(2)},\dots] [[\zeta]]$ for all~$i,j\geq 0$.

\begin{lem}\label{coef}
The following formulae hold true:
\begin{align}
& D(\lambda)(u^{(i)})
=\frac{1}{\lambda}\left(\p_x^iB_{0,0}+\sum_{j=1}^i\binom{i}{j}B_{j-1,i-j+1}\right), \quad \forall\,i\geq 0, \label{c111}\\
& \e^{-2}\left[J_2(\lambda)(u^{(i)})J_2(\lambda)(u^{(j)})\right]_{{\rm reg},-}=\frac{1}{\lambda}B_{i+1,j+1}, \quad \forall\, i,j\geq 0. \label{c222}
\end{align}
\end{lem}
\begin{prf}
Using Lemma~\ref{prp-1} and the definition of~$J_2(\lambda)$ given in Lemma~\ref{virasoro-total}, we have
\beq\label{J2pxf0}
\e^{-1}J_2(\lambda)\left(\p_x\mathcal{F}_0\right)+\frac{1}{\sqrt{\lambda}} =
\frac{1}{\sqrt{\lambda}}\left(E_0+\sqrt{\frac{K_2}{h}}\sum_{\alpha=1}^{K_1-1}E_{\alpha}
+\sqrt{\frac{K_1}{h}}\sum_{\alpha=-(K_2-1)}^{-1}E_{\alpha}\right).
\eeq
Acting~$\p_x^2$ on both sides of the equation~\eqref{symmetry-0} yields
\[D(\lambda)(u)=\left[\left(\e^{-1}J_2(\lambda)(\p_x\mathcal{F}_0)+\frac{1}{\sqrt{\lambda}}\right)^2\right]_{{\rm reg},-}
=\frac{1}{\lambda} B_{0,0}, \]
which gives the validity of~\eqref{c111} for $i=0$. Assume that the formula~\eqref{c111} is true for $i=k$; 
then for $i=k+1$ we have
\begin{align*}
 &D(\lambda)(u^{(k+1)})
= \p_x D(\lambda)(u^{(k)})+[D(\lambda),\,\p_x](u^{(k)})\\
=& \frac{1}{\lambda}\left(\p_x^{k+1}B_{0,0}+\sum_{j=1}^k\binom{k}{j}\left(B_{j,k-j+1}+B_{j-1,k-j+2}\right)\right) \nn\\
&+ \e^{-1}\Bigl[ J_2(\lambda)(\p_x\mathcal{F}_0)
+\tfrac{1}{\sqrt{\lambda}} J_2(\lambda)(u^{(k)})\Bigr]_{{\rm reg},-}\\
=& \frac{1}{\lambda}\left(\p_x^{k+1}B_{0,0}+\sum_{j=1}^{k+1}\binom{k+1}{j}B_{j-1,k-j+2}\right).
\end{align*}
Hence by induction we arrive at the formula~\eqref{c111}. 
The formula~\eqref{c222} can be proved by using  \eqref{zh-2} and \eqref{J2pxf0}. The lemma is proved.
\end{prf}

Let us now introduce a family of polynomials $f_{i,j} \in \C\bigl[u^{(1)},u^{(2)},\dots\bigr]$, $i,j\geq 0$ by 
requiring the validity of the identity
\beq\label{deff}
\p_x^i y(\zeta)=\sum_{j=0}^i f_{i,j}\left(\zeta\p_\zeta\right)^j y(\zeta)
\eeq
for any smooth function $y(\zeta)$ of~$\zeta=e^{u(x,{\bf s})}/\lambda$. Clearly, $f_{i,j}$ vanishes if $i<j$. 
For $i\geq j$, the polynomials $f_{i,j}$ can be uniquely determined by the following recursive relations:
\begin{align*}
&f_{i,0}=\delta_{i,0},\\
&f_{i+1,j+1}=\p_x f_{i,j+1}+u^{(1)} f_{i,j}.
\end{align*}
Explicit expressions for~$f_{i,j}$ will be given in Section~\ref{section4}.
The functions $B_{i, j}$ defined in~\eqref{zh-2} can now be written as
\begin{equation}\label{f-0} 
B_{i,j}=\sum_{k=0}^i\sum_{l=0}^jf_{i,k}f_{j,l}\widetilde{B}_{k,l},
\end{equation}
where
\begin{equation}
\widetilde{B}_{k,l}:= \sum_{\alpha,\beta\in I} \left(\zeta\p_{\zeta}\right)^k \left(E_{\alpha}\right) 
G^{\alpha\beta} \left(\zeta\p_{\zeta}\right)^l \left(E_{\beta}\right), \label{widetb}
\end{equation}
and $E_\alpha$ are defined by~\eqref{eeplus}.

\begin{lem}\label{algo}
The functions $\widetilde{B}_{i,j}$, $i,j\geq 0$ defined in~\eqref{widetb} satisfy the relations
\begin{align}
& \sum_{k\geq0}z^{-k}\widetilde{B}_{0,k} =\sqrt{z}\Phi(z)\frac{1}{1-K \zeta e^{-\p_z}}\left(\frac{1}{\sqrt{z}\Phi(z)}\right), \label{f-1} \\
& \widetilde{B}_{i,j}=\widetilde{B}_{j,i},\quad
\zeta\p_\zeta \widetilde{B}_{i,j}=\widetilde{B}_{i+1,j}+\widetilde{B}_{i,j+1}. \label{recur-f} 
\end{align}
Moreover, $\widetilde{B}_{i,j}$ are polynomials of~$(1-K\zeta)^{-1}$ with degrees less than or equal to~$i+j+1$.
\end{lem}
\begin{prf} 
By substituting~\eqref{eeplus} into~\eqref{widetb} we obtain
\begin{equation}
\widetilde{B}_{0,k}=\sum_{n\geq0}A_{k,n}\zeta^n, 
\end{equation}
where
\begin{align*}
A_{k,n}:=&c_{h n} \delta_{k,0} + \sum_{\ell=1}^{n-1}b_{h\ell}^k c_{h\ell} c_{h(n-\ell)} + \frac{K_2}{h}\sum_{\alpha=1}^{K_1-1}\sum_{\ell=0}^{n-1}b_{\alpha+h\ell}^k c_{\alpha+h\ell}c_{K_1-\alpha+h(n-1-\ell)} \\
& + \frac{K_1}{h}\sum_{\alpha=-(K_2-1)}^{-1}\sum_{\ell=0}^{n-1}b_{\alpha+h\ell}^k c_{\alpha+h\ell}c_{-\alpha-K_2+h(n-1-\ell)}.\\
\end{align*}
Using Lemma~\ref{const} we find that
\begin{align}
K^{-n} A_{k,n}
=&\sum_{\alpha\in I}\sum_{\ell=0}^{n-1}b_{\alpha+h\ell}^{k-1}\Res_{z=b_{\alpha+h\ell}}V_n(z)
 +\frac{\Res_{z=n}V_n(z)}{n}\delta_{k,0}.
\end{align}
In particular, we have 
\[A_{0,n}=K^n\left(1-V_n(0)+\frac1n\Res_{z=n}V_n(z)\right)=K^n.\]
Hence
\begin{align}
&\sum_{k\geq0}z^{-k}\widetilde{B}_{0,k}  =  \sum_{k,n\geq0}z^{-k}A_{k,n}\zeta^n  = \sum_{n\geq0}V_n(z) \left(K\zeta\right)^n \nn\\
= & \sum_{n\geq 0} \frac{\Phi(z)}{\Phi(z-n)}\sqrt{\frac{1}{1-n/z}} \left(K\zeta\right)^n = \sqrt{z}\Phi(z)\frac{1}{1-K\zeta e^{-\p_z}}\left(\frac{1}{\sqrt{z}\Phi(z)}\right). \label{324}
\end{align}
Note that 
\beq
\frac1{1-K \zeta e^{-\p_z}} = \sum_{n\geq 0} \sum_{k=1}^{n+1} \frac{Q(n,k)}{(1-K\zeta)^k} \frac{\p_z^n}{n!},\quad Q(n,k)=\frac{1}{k}\sum_{i=1}^k(-1)^{i-1}\binom{k}{i}i^{n+1}. \label{QQQQ}
\eeq
So from~\eqref{recur-f} and~\eqref{324} it follows that $\widetilde{B}_{i,j}$ belong to $\C\bigl[(1-K\zeta)^{-1}\bigr]$ with 
degrees less than or equal to~$i+j+1$. 
The lemma is proved.
\end{prf}

\begin{rmk}
Introduce the trinomial curve 
\beq \label{curve}
X Y^{K_1+K_2} - Y^{K_2} + 1 = 0 . 
\eeq
Near $X=0$ we have two Puiseux series solutions:
\begin{align}
& Y_- = \sum_{m\geq 0} c_{-,m} X^m , \quad c_{-,0}=1, \nn \\
& Y_+ = X^{-\frac1{K_1}} \sum_{m\geq 0} c_{+,m} X^{m\frac{K_2}{K_1}},  \quad c_{+,0}=1. \nn
\end{align}
By using the Lagrange inversion formula we obtain
\[
(Y_-)^{K_2} = \sum_{m\geq 0} \binom{m \frac{K_1+K_2}{K_2}}{m} \frac{X^m}{1+m\frac{K_1}{K_2}} , 
\quad (Y_+)^{-K_1} = \sum_{m\geq 0} \binom{m \frac{K_1+K_2}{K_1}}{m} \frac{X^{1+m \frac{K_2}{K_1}}}{1+m\frac{K_2}{K_1}}. 
\]
We also observe the validity of the following relations:
\begin{align}
& \sum_{\alpha=0}^{K_1-1} E_{\alpha} = \frac{d}{dX} \left((Y_+)^{-K_1}\right), 
\quad  \sum_{\alpha=-(K_2-1)}^{0} E_{\alpha} = (Y_-)^{K_2} + \frac{K_1}{K_2} X \frac{d}{dX} \left((Y_-)^{K_2}\right) 
\end{align}
with $X=\zeta^{1/K_2}$.
We conjecture that the curve~\eqref{curve} can be used 
to give an alternative way of computing the special cubic Hodge integrals
by using the topological recursion of Chekhov--Eynard--Orantin type. 
We note that the $K_2=1$ case was considered by Bouchard--Klemm--Mari\~no--Pasquetti~\cite{BKMP}. 
\end{rmk}
\begin{cor}\label{J2-F0}
The following formula holds true:
\beq
\e^{-2} \Bigl[J_2(\lambda)^2 \left(\mathcal{F}_0\right)\Bigr]_{{\rm reg},-}
=\frac{1}{\lambda}
\left(\frac{1}{8(1-K\zeta)^2}-\left(\frac{1}{8}-\frac{\sigma_1}{12}\right)\frac{1}{1-K \zeta}-\frac{\sigma_1}{12}\right).\label{zh-3}
\eeq
\end{cor}

\begin{prf} By using the definition of~$J_2(\lambda)$ given in Lemma~\ref{virasoro-total} and by using Lemma~\ref{prp-1} we obtain
\beq\label{zh-4}
\e^{-2} \left[J_2(\lambda)^2(\mathcal{F}_0)\right]_{{\rm reg},-}=\int \frac{1}{\lambda}\frac{B_{1,1}}{u^{(1)}} dx = 
\frac{1}{\lambda}\int \frac{\widetilde{B}_{1,1}}{K\zeta}d(K\zeta).
\eeq
Here we recall that $B_{1,1}$ is defined in~\eqref{zh-2}, and
the integration constant is chosen such that the left hand side tends to zero as $\lambda\rightarrow \infty$. 
From Lemma~\ref{algo} it follows that
\begin{align*}
&\widetilde{B}_{1,1}
=\frac{1}{4}\frac{1}{(1-K\zeta)^3}
-\left(\frac{3}{8}-\frac{\sigma_1}{12}\right)\frac{1}{(1-K \zeta)^2}
+\left(\frac{1}{8}-\frac{\sigma_1}{12}\right)\frac{1}{1-K \zeta},
\end{align*}
which, together with~\eqref{zh-4}, leads to the formula~\eqref{zh-3}. The corollary is proved.\
\end{prf}

By using Lemmas~\ref{lemma35}, \ref{coef} and Corollary~\ref{J2-F0} 
we arrive at the following theorem.
\begin{thm}\label{main} 
The series $\Delta \mathcal{F}$ given in \eqref{zh-5}, \eqref{zh-6} satisfies the Dubrovin--Zhang loop equation
\begin{align} 
& \sum_{i\geq0}\frac{\p \Delta \mathcal{F}}{\p u^{(i)}} \p_x^i \Theta 
+\sum_{i\geq1} \sum_{j=1}^i \sum_{\alpha,\beta\in I} \frac{\p \Delta \mathcal{F}}{\p u^{(i)}}\binom{i}{j} 
\p_x^{j-1}(E_{\alpha})G^{\alpha\beta}\p_x^{i-j+1}(E_{\beta}) 
 \nn \\
& =  \frac{\Theta^2}{16}-\left(\frac{1}{16}-\frac{\sigma_1}{24}\right)\Theta + \e^2\sum_{i\geq0} 
\frac{\p \Delta \mathcal{F}}{\p u^{(i)}} \p_x^{i+2}\left(\frac{\Theta^2}{16}- \left(\frac{1}{16}-\frac{\sigma_1}{24}\right)\Theta\right) \nn \\
& \quad +\frac{\e^2}{2}\sum_{i,j\geq0} \sum_{\alpha,\beta \in I}  \left(\frac{\p^2\Delta \mathcal{F}}{\p u^{(i)} \p u^{(j)}}
+\frac{\p \Delta \mathcal{F}}{\p u^{(i)}}\frac{\p \Delta \mathcal{F}}{\p u^{(j)}}\right) 
\p_x^{i+1}(E_{\alpha})G^{\alpha\beta}\p_x^{j+1}(E_{\beta}), \label{loop-eq}
\end{align}
where $\Theta= \frac1{1- K \zeta}$ with $K=h^h K_1^{-K_1} K_2^{-K_2}$, and 
$\sigma_1= -(p+q+r)=\frac1{h}-\frac1{K_1}-\frac1{K_2}$.
\end{thm}

Note that each side of the loop equation~\eqref{loop-eq} is a power series of~$\Theta$. 
It is understood that this equation for~$\Delta \mathcal{F}$ holds true identically in~$\Theta$. 
It is easy to check that $\mathcal{F}_g$ also satisfy the equations
\beq
\mathcal{F}_1= \frac1{24} \log u^{(1)} + \frac{\sigma_1}{24} u^{(0)}, \quad 
\sum_{j=1}^{3g-2} j u^{(j)} \frac{\p \mathcal{F}_g}{\p u^{(j)}} = (2g-2) \mathcal{F}_g, \quad g\geq 2.  \label{eqfix} 
\eeq

\begin{prp}\label{uniquenessprop}
The solution to equations \eqref{loop-eq}--\eqref{eqfix} is unique. 
\end{prp}
\begin{prf}
Expanding~\eqref{loop-eq} with respect to~$\e$ and comparing
 the coefficients of powers of~$\e$ we find that~\eqref{loop-eq} is equivalent to the following 
 equations:
\begin{align*}
&\sum_{i\geq 0}\left(\p_x^i B_{0,0}
+\sum_{j=1}^i\binom{i}{j}B_{j-1,i-j+1}\right)
\frac{\p \mathcal{F}_1}{\p u^{(i)}}=\frac{\Theta^2}{16}-\left(\frac{1}{16}-\frac{\sigma_1}{24}\right) \Theta , \nn\\
&\sum_{i\geq 0}\left(\p_x^i B_{0,0}
+\sum_{j=1}^i\binom{i}{j}B_{j-1,i-j+1}\right)\frac{\p \mathcal{F}_g}{\p u^{(i)}}=\sum_{i\geq0}\p_x^{i+2}
\left(\frac{\Theta^2}{16}
- \left(\frac{1}{16}-\frac{\sigma_1}{24}\right)\Theta\right)
\frac{\p\mathcal{F}_{g-1}}{\p u^{(i)}}\\
& \quad + \frac12\sum_{i,j\geq0}B_{i+1,j+1}
\left(\frac{\p^2 \mathcal{F}_{g-1}}{\p u^{(i)}\p u^{(j)}}
+ \sum_{k=1}^{g-1}\frac{\p \mathcal{F}_k}{\p u^{(i)}}\frac{\p \mathcal{F}_{g-k}}{\p u^{(j)}}\right), \quad g\geq 2. \nn
\end{align*}
Here we used the formulae~\eqref{zh-2}.
By using the fact that 
$\p \mathcal{F}_g/\p u^{(i)}=0$ for $i\geq 3g-1$, 
and that $B_{i,j}$ 
are polynomials in~$\Theta$ of degrees~$i+j+1$, we arrive at the following system of equations:
\begin{align*}
&\left(\Theta,\cdots,\Theta^{3g-1}\right) M_g 
\left(\frac{\p \mathcal{F}_g}{\p u^{(0)}},\cdots,\frac{\p \mathcal{F}_g}{\p u^{(3g-2)}}\right)^T=\left(\Theta,\cdots,\Theta^{3g-1}\right) N_g, \quad g\geq 1,
\end{align*}
where $M_g$ is an invertible upper triangular $(3g-1)$ by $(3g-1)$ matrix, and $N_g$ 
is a column vector. All the elements of the matrix and the vector are polynomials of ~$1/u^{(1)}, u^{(0)}, u^{(1)}$, $\cdots$. 
So the gradient of $\mathcal{F}_g$ is uniquely determined. 
The proposition then follows from \eqref{eqfix}.
\end{prf}

\section{Loop equation: the general case}\label{section4}
In this section we drop the rational condition~\eqref{rational-condition} and consider the general case 
when $p, q, r$ are arbitrary complex numbers satisfying the {\it local Calabi--Yau condition \eqref{calabi-yau}}. 
We will give and prove  
the general version of the loop equation for the 
corresponding Hodge free energies.

As it is pointed out in Lemma~\ref{jetvariablesrep}, there are polynomials 
$H_g(z_0, z_1, \dots, z_{3g-2};p,q,r)$ in $p, q, r$ satisfying the relations~\eqref{HHgHg}
and the relations~\eqref{eq2233} (see \cite{DLYZ-1}), where $g\geq1$. It is then clear that 
$H_g(z_0, \dots, z_{3g-2};p,q,r)$, $g\geq1$ are also polynomials of $\sigma_1,\sigma_3$. For simplicity, we will denote them by
$H_g(z_0,\dots,z_{3g-2}; \sigma_1,\sigma_3)$.

We will need the following lemma.
\begin{lem}\label{elem}
Let $\mathcal{P}(\sigma_1,\sigma_3)$ be a polynomial in $\C[\sigma_1,\sigma_3]$. 
If $\mathcal{P}$ vanishes for the values 
\[(\sigma_1,\sigma_3)=\left(\frac1{K_1+K_2}-\frac1{K_1}-\frac1{K_2}, 
\frac2{(K_1+K_2)^3}- \frac{2}{K_1^3} - \frac{2}{K_2^3}\right),\] 
where 
$K_1,K_2 \in \N$, $(K_1,K_2)=1$, then $\mathcal{P}\equiv 0$. 
\end{lem}
\begin{prf} 
Suppose $\mathcal{P}\not\equiv 0$. 
Fix $K_1$ to be any positive integer. We observe that the points 
\[\left(\frac1{K_1+K_2}-\frac1{K_1}-\frac1{K_2}, \frac2{(K_1+K_2)^3}- \frac2{K_1^3} - \frac{2}{K_2^3}\right), \quad K_2\in \N\] 
belong to the irreducible algebraic curve 
\[2 K_1^3 \sigma_1^3 - 6 K_1 \sigma_1 - 6 - K_1^3 \sigma_3 =0\]
on the $(\sigma_1,\sigma_3)$ plane. It is easy to see that there are infinitely 
many of such points with $(K_1,K_2)=1$. Hence the polynomial 
\[2 K_1^3 \sigma_1^3 - 6 K_1 \sigma_1 - 6 - K_1^3 \sigma_3\]
must divide $\mathcal{P}(\sigma_1,\sigma_3)$.
This contradicts with the fact that a polynomial has a finite degree. The lemma is proved.
\end{prf}

It follows from Lemma~\ref{elem} that if a polynomial~$P$ 
in $\C[\sigma_1,\sigma_3][z_2,z_3,\dots]\bigl[z_1,z_1^{-1}\bigr]$
is equal to~$H_g(z_0,\dots,z_{3g-2};\sigma_1,\sigma_3)$ 
for 
\[(\sigma_1,\sigma_3)=\left(\frac1{K_1+K_2}-\frac1{K_1}-\frac1{K_2}, 
\frac2{(K_1+K_2)^3}- \frac2{K_1^3} - \frac{2}{K_2^3}\right)\] 
with 
$K_1,K_2$ being arbitrary coprime positive integers, 
then $P\equiv H_g(z_0,\dots,z_{3g-2};\sigma_1,\sigma_3)$, where $g\geq 2$.
For $g=1$, we have the explicit expression
\[H_1(z_0,z_1;\sigma_1,\sigma_3)= \frac1{24} \log z_1 + \frac{\sigma_1}{24}z_0.\]


To prove Theorem~\ref{refineloopmain}, we futher introduce some notations.
Let $\mathcal{B}_{n,k}(X_1,\dots,X_{n-k+1})$ be the exponential Bell polynomials. They can be defined via the generating function
$$
\sum_{n,k\geq0}\mathcal{B}_{n,k} \left(X_1,\dots,X_{n-k+1}\right)\frac{y^n}{n!}z^k=\exp\left(z\sum_{j\geq1} X_j\frac{y^j}{j!}\right).
$$
Let $\mathcal{B}_n$ denote the complete Bell polynomials, i.e.
\[
\mathcal{B}_0:=1,\quad \mathcal{B}_n(X_1,\dots,X_n) := \sum_{k=1}^n \mathcal{B}_{n,k} \left(X_1,\dots,X_{n-k+1}\right), \quad n\geq 1.
\]
Now we define $\widetilde{P}_{0,n}(\xi;\sigma_1,\sigma_3) \in \C[\sigma_1,\sigma_3] \bigl[(1-\xi)^{-1}\bigr]$, $n\geq0$ by the equation
\begin{equation}\label{zh-7}
\sum_{n\geq0}z^{-n}\widetilde{P}_{0,n}(\xi; \sigma_1,\sigma_3)
=\sqrt{z}\Phi(z;\sigma_1,\sigma_3)\frac{1}{1-\xi e^{-\p_z}}\left(\frac{1}{\sqrt{z}\Phi(z;\sigma_1,\sigma_3)}\right),
\end{equation}
where 
\[\Phi(z;\sigma_1,\sigma_3)=\exp\left(-\sum_{i\geq1}\frac{B_{2i}}{2i(2i-1)}(p^{2i-1}+q^{2i-1}+r^{2i-1})z^{1-2i}\right).\]
Denote $\Phi=\Phi(z;\sigma_1,\sigma_3)$. Then the equation~\eqref{zh-7} can be written more explicitly as follows:
\begin{align}
& \sum_{n\geq0}z^{-n}\widetilde{P}_{0,n}(\xi;\sigma_1,\sigma_3)=\sum_{n\geq0}\sum_{k=1}^{n+1}\frac{Q(n,k)}{(1-\xi)^k}
\sum_{m=0}^n \frac{(-1)^{m}(2m-1)!!}{2^{m}m! (n-m)!}z^{-m}
\Phi\, \p_z^{n-m}\left(\frac{1}{\Phi}\right), \label{B0xi} 
\end{align}
where $Q(n,k)$ are the numbers defined in~\eqref{QQQQ}, and 
\begin{align}
&\Phi \, \p_z^\ell\left(\frac{1}{\Phi}\right)
=\mathcal{B}_\ell \left(-\p_z\log\Phi,\dots,-\p_z^\ell\log\Phi\right), \quad \ell\geq0.\label{phiphi}
\end{align}
We define, for $i,j\geq0$, $\widetilde{P}_{i,j}=\widetilde{P}_{i,j}(\xi;\sigma_1,\sigma_3)$ by the following recursion:
\begin{align*}
\widetilde{P}_{i,j} &= \widetilde{P}_{j,i}, \\
\xi\p_{\xi}\widetilde{P}_{i,j} &= \widetilde{P}_{i+1,j}+\widetilde{P}_{i,j+1}.
\end{align*}

Note that the functions $f_{i,k}$ defined in \eqref{deff} are just the Bell polynomials
$\mathcal{B}_{i,k}(z_1,\dots,z_{i-k+1})$. By using these functions we 
introduce the following notations:
$$
P_{i,j}=\sum_{k=0}^i\sum_{\ell=0}^j \mathcal{B}_{i,k}(z_1,\dots,z_{i-k+1}) \mathcal{B}_{j,\ell}(z_1,\dots,z_{j-\ell+1})\widetilde{P}_{k,\ell}(\xi;\sigma_1,\sigma_3),\quad i,j\geq 0.
$$

We are now ready to prove Theorem~\ref{refineloopmain}. 

\smallskip

\begin{prfn}{Theorem~\ref{refineloopmain}}
It is easy to verify that the functions $P_{i,j}$ defined above 
coincide with the functions $B_{i,j}$ that are given in~\eqref{f-0} when $\xi =K \zeta$ and
\begin{align*}
& (\sigma_1, \sigma_3) = \left(\frac1{h}-\frac1{K_1}-\frac1{K_2}, 
\frac2{h^3}- \frac{2}{K_1^3} - \frac{2}{K_2^3}\right),\\
& (z_1,z_2,\dots) = \left(u^{(1)},u^{(2)},\dots\right),
\end{align*}
where $K_1,K_2$ are arbitrary coprime positive integers. 
It then follows from Theorem~\ref{main} that 
$\Delta H:=\sum_{g\geq 1} \e^{2g-2} H_g(z_0,\dots,z_{3g-2};\sigma_1,\sigma_3)$ satisfies 
equation~\eqref{loop-eq-0} for the particular values of $(\sigma_1,\sigma_3)$, where 
$H_g(z_0,\dots,z_{3g-2};\sigma_1,\sigma_3)$ are introduced in the beginning of this section.
By definition we also see that $P_{i,j}\in \C[\sigma_1,\sigma_3][(1-\xi)^{-1};z_1,z_2,\dots]$.
Then by using Lemma~\ref{elem} we find that $\Delta H$ 
is a solution to the loop equation~\eqref{loop-eq-0}. 
As we do in the proof of Proposition~\ref{uniquenessprop}, we deduce from~\eqref{loop-eq-0} 
the following equation for each $g\geq 1$:
\begin{align*}
&\left(\Theta,\cdots,\Theta^{3g-1}\right) \widetilde{M}_g(\sigma_1,\sigma_3) 
\left(\frac{\p H_g}{\p z_0},\cdots,\frac{\p H_g}{\p z_{3g-2}}\right)^T
=\left(\Theta,\cdots,\Theta^{3g-1}\right) \widetilde{N}_g(\sigma_1,\sigma_3),
\end{align*}
where $\widetilde{M}_g(\sigma_1,\sigma_3)$ is an invertible $(3g-1)$ by $(3g-1)$ matrix, 
$\widetilde{N}_g(\sigma_1,\sigma_3)$ is 
a column vector, and their entries are polynomials 
of $\sigma_1,\sigma_3$. Therefore we verified that the gradient of~$H_g$ is uniquely determined by~\eqref{loop-eq-0}.
The theorem is proved.  
\end{prfn}

\begin{emp}
By taking the coefficient of $\e^0$ in \eqref{loop-eq-0}, we obtain the following equation for~$H_1$:
$$
P_{0,0}\frac{\p H_1}{\p z_0}
+\left(\p P_{0,0}+P_{0,1}\right)\frac{\p H_1}{\p z_1}
=\frac{\Theta^2}{16}
-\left(\frac{1}{16}-\frac{\sigma_1}{24}\right) \Theta.
$$
The coefficients read
\[P_{0,0}=\Theta,\quad P_{0,1}
=\frac{z_1}{2}\left(\Theta^2-\Theta\right).\]
So we have
\begin{align*}
&\frac{3 z_1}{2}\frac{\p H_1}{\p z_1} \Theta^2
+\left(\frac{\p H_1}{\p z_0}-\frac{3 z_1}{2}\frac{\p H_1}{\p z_1}\right)
\Theta=\frac{\Theta^2}{16}
-\left(\frac{1}{16}-\frac{\sigma_1}{24}\right)\Theta,
\end{align*}
which gives
\[\frac{\p H_1}{\p z_0}=\frac{\sigma_1}{24},\quad
\frac{\p H_1}{\p z_1}=\frac{1}{24 z_1}.\]
Hence 
we obtain~\eqref{explicitH1}.
In a similar way, we obtain formula~\eqref{explicitH2} and the following expression of~$H_3$:
\begin{align}
H_3
=&\frac{1}{82944}\frac{z_7}{z_1^3}-\frac{7}{46080}\frac{z_2 z_6}{z_1^4}
 -\frac{53}{161280}\frac{z_3 z_5}{z_1^4}
 +\frac{353}{322560}\frac{z_2^2 z_5}{z_1^5}
 -\frac{103}{483840}\frac{z_4^2}{z_1^4}\nn\\
 &+\frac{1273}{322560 }\frac{ z_2 z_3 z_4}{z_1^5}
-\frac{83 }{15120 }\frac{z_2^3 z_4}{z_1^6}+\frac{59 }{64512 }\frac{z_3^3}{z_1^5}
 -\frac{83  }{7168 }\frac{z_2^2 z_3^2}{z_1^6}+\frac{59}{3024 }\frac{ z_2^4 z_3}{z_1^7}\nn\\
&-\frac{5 }{648 }\frac{z_2^6}{z_1^8} +\frac{7 \sigma_1 }{138240}\frac{z_6}{ z_1^2}
 -\frac{383 \sigma_1 }{967680 }\frac{z_2 z_5}{z_1^3}
+\frac{41 \sigma_1^2 }{580608 }\frac{z_5}{z_1}-\frac{689 \sigma_1 }{967680}\frac{z_3 z_4}{ z_1^3}
+\frac{185 \sigma_1 }{96768 }\frac{z_2^2 z_4}{z_1^4}\nn\\
&-\frac{373 \sigma_1^2 }{1451520 }\frac{z_2 z_4}{z_1^2}
 +\left(\frac{23 \sigma_1^3}{580608}-\frac{11 \sigma_3}{2903040}\right)z_4 
+\frac{869 \sigma_1 }{322560 }\frac{z_2 z_3^2}{z_1^4}
-\frac{61 \sigma_1^2 }{322560 }\frac{z_3^2}{z_1^2}\nn\\
& -\frac{9343 }{1451520 }\frac{z_2^3 z_3}{z_1^5}+\frac{151 \sigma_1^2 }{207360 }\frac{z_2^2 z_3}{z_1^3}
 -\left(\frac{19 \sigma_1^3}{1451520}-\frac{19 \sigma_3}{2903040}\right)\frac{z_2 z_3}{z_1} +\frac{131 \sigma_1 }{45360 }\frac{z_2^5}{z_1^6}\nn\\
&-\frac{19 \sigma_1^2 }{53760 }\frac{z_2^4}{z_1^4} +\left(\frac{41 \sigma_1^4}{4354560}-\frac{41 \sigma_1 \sigma_3}{8709120}\right)z_1 z_3
 +\left(\frac{\sigma_1^3}{108864}-\frac{\sigma_3}{217728}\right)\frac{z_2^3}{z_1^2} \nn \\
& +\left(\frac{31 \sigma_1^4}{4354560}-\frac{31 \sigma_1 \sigma_3}{8709120}\right)z_2^2
 -\left(\frac{\sigma_1^6}{13063680}-\frac{\sigma_1^3 \sigma_3}{13063680}
 +\frac{\sigma_3^2}{52254720}\right)z_1^4 . \label{explicitH3}
 \end{align}
\end{emp}

\medskip

\noindent E-mails: \\
liusq@tsinghua.edu.cn\\
diyang@ustc.edu.cn\\
youjin@tsinghua.edu.cn\\
zhouch@ustc.edu.cn


\begin{thebibliography}{99}
\bibitem{BKMP}
Bouchard, V., Klemm, A.,  Mari\~no, M., Pasquetti, S., 
Remodeling the B-model. Comm. Math. Phys. {\bf 287} (2009), 117--178.

\bibitem{BPS1}
Buryak, A., Posthuma, H., Shadrin, S., 
A polynomial bracket for the Dubrovin-Zhang hierarchies. J. Differential Geom. {\bf 92} (2012), 153--185.

\bibitem{BPS2}
Buryak, A., Posthuma, H., Shadrin, S., 
 On deformations of quasi-Miura transformations and the Dubrovin-Zhang bracket. J. Geom. Phys. {\bf 62} (2012), 1639--1651.
 
\bibitem{DVV}
Dijkgraaf, R., Verlinde, H., Verlinde, E., Loop equations and Virasoro constraints in 
nonperturbative two-dimensional quantum gravity. Nucl. Phys. B, {\bf 348} (1991), 435-456.

\bibitem{Du1}
Dubrovin, B., Geometry of 2D topological field theories. 
In ``Integrable Systems and Quantum Groups" (Montecatini Terme, 1993). 
Editors: Francaviglia, M., Greco, S. 
Springer Lecture Notes in Math.~\textbf{1620} (1996), 120--348.

\bibitem{Du2}
Dubrovin, B., On almost duality for Frobenius manifolds. 
Translations of the American Mathematical Society-Series 2, {\bf 212} (2004), 75--132.

\bibitem{DLYZ-1} Dubrovin, B., Liu, S.-Q., Yang, D., Zhang, Y., 
 Hodge integrals and tau-symmetric integrable hierarchies of Hamiltonian evolutionary PDEs. 
 Adv. Math. {\bf 293} (2016), 382--435.
 

 \bibitem{DLYZ-2}
Dubrovin, B., Liu, S.-Q., Yang, D., Zhang, Y., Hodge-GUE correspondence and the discrete KdV equation. 
Comm. Math. Phys. {\bf 379} (2020), 461--490.

\bibitem{DY} 
Dubrovin, B., Yang, D.,  On cubic Hodge integrals and random matrices.
Commun. Number Theory Phys. {\bf 11} (2017), 311--336.

\bibitem{DZ}
Dubrovin, B., Zhang, Y., Normal forms of hierarchies of integrable PDEs, Frobenius manifolds and Gromov-Witten invariants. arXiv:math/0108160.

\bibitem{FP}
Faber, C., Pandharipande, R., Hodge integrals and Gromov-Witten theory. Invent. Math. {\bf 139} (2000), 173--199.

\bibitem{G}
Givental, A.B., Gromov-Witten invariants and quantization of quadratic Hamiltonians.  
Mosc. Math. J. {\bf 1} (2001), 551--568.

\bibitem{GV}
Gopakumar, R., Vafa, C., On the gauge theory/geometry correspondence. Adv. Theor. Math. Phys. {\bf 5} (1999), 1415--1443.

\bibitem{GP} 
Graber, T., Pandharipande, R., Localization of virtual classes. Invent. Math. {\bf 135} (1999), 487--518.

\bibitem{HKQ} 
Huang, M.-X., Klemm, A., Quackenbush, S., 
Topological string theory on compact Calabi-Yau: modularity and boundary conditions. 
In ``Homological mirror symmetry", 45--102. Springer, Berlin, Heidelberg, 2009.

\bibitem{Kontsevich} 
Kontsevich, M., Intersection theory on the moduli space of curves and the matrix Airy function. 
Comm. Math. Phys. {\bf 147} (1992), 1--23.

\bibitem{LLZ}
Liu, C.-C. M., Liu, K., Zhou, J., A proof of a conjecture of Mari\~no--Vafa on Hodge integrals. 
J. Differential Geom. {\bf 65} (2003), 289--340.

\bibitem{LYZZ}
Liu, S.-Q., Yang, D., Zhang, Y., Zhou, C., The Hodge-FVH correspondence. 
J. Reine Angew. Math., to appear. Preprint arXiv:1906.06860.

\bibitem{LZZ}
Liu, S.-Q., Zhang, Y., Zhou, C., Fractional Volterra Hierarchy. Lett. Math. Phys. {\bf 108} (2018), 261--283.

\bibitem{MV}
Mari\~no, M., Vafa, C., Framed knots at large N. Orbifolds in Mathematics and Physics (Madison,
WI, 2001), Contemp. Math. {\bf 310} (2002), 185--204.

\bibitem{M} 
Mumford, D., Towards an enumerative geometry of the moduli space of curves. In ``Arithmetic and geometry"~II, 271--328,  Progr. Math., 36, Birkh\"auser, Boston, MA, 1983.

\bibitem{OP}
Okounkov, A., Pandharipande, R., Hodge integrals and invariants of the unknot. 
Geom. Topol. {\bf 8} (2004), 675--699.

\bibitem{WW}
Whittaker, E.T., Watson, G.N., A Course of Modern Analysis, 4th edn. Cambridge University Press, Cambridge, 1963.

\bibitem{Witten} Witten, E., Two-dimensional gravity and intersection theory on moduli space. Surveys in differential geometry (1991), 243--320. Lehigh Univ., Bethlehem, PA.

\bibitem{Z}
Zhou, J., On recursion relation for Hodge integrals from cut-and-join equations. Preprint, 2009.

\end{thebibliography}
\end{document}